 \documentclass[11pt]{amsart}  
\usepackage{amssymb,color} 
\newcommand{\eq}[1]{{\begin{equation}#1\end{equation}}}
\newcommand{\spl}[1]{{\begin{split}#1\end{split}}}
\newcommand{\eqn}[1]{\begin{equation*}#1\end{equation*}}

\newtheorem{thm}{Theorem}[section]
\newtheorem{cor}[thm]{Corollary}
\newtheorem{prop}[thm]{Proposition}
\newtheorem{lem}[thm]{Lemma}

\newenvironment{prf}[1]
   {{\noindent \bf Proof of {#1}.}}{\hfill \qed}

\renewcommand{\r}{\rho}

\newcommand{\re}{\mathbb{R}}

\newcommand{\A}{\mathbb{A}} 
\newcommand{\B}{\mathbb{B}} 
     
\newcommand{\pt}{\partial}   
\newcommand{\al}{\alpha}

\newcommand{\gm}{\gamma}
\newcommand{\ep}{\varepsilon}
\newcommand{\lam}{\lambda}

\newcommand{\del}{\delta}

\renewcommand{\(}{|\!|\!| }
\renewcommand{\)}{|\!|\!|}
\def\<{\langle}
\renewcommand{\>}{\rangle}

\begin{document}
\title{Global well-posedness for one dimensional Chern-Simons-Dirac 
system in $L^p$}
\date{2012,9/27}
\maketitle  
\vskip5mm
{\small
\vskip3mm\hskip7mm
{\normalsize \sf Shuji Machihara}
\vskip1mm\hskip1.5cm
Faculty of Education, School of Mathematics, \\
\hskip1.5cm\indent
Saitama University \\ 
\hskip1.5cm\indent
Saitama 338-8570, JAPAN
\vskip3mm\hskip7mm
{\normalsize \sf Takayoshi Ogawa}
\vskip1mm\hskip1.5cm
Mathematical Institute,
Tohoku University \\ 
\hskip1.5cm\indent
Sendai 980-8578, JAPAN
}

\baselineskip=5.6mm

\vskip5mm

\baselineskip 5.5mm

\begin{abstract}
Time global wellposedness in $L^{p}$ for the Chern-Simons-Dirac equation 
in $1+1$ dimension is discussed. 
The two types of quadratic terms are treated, null case and non-null case.
The standard iteration arguments, 
different settings correspond to each cases respectively,
are used for the proof.
For the critical case in $L^{1}$, 
the mass concentration phenomena of the solutions is denied to show
the time global solvability. The intrinsic estimate plays an important
role in the proof. These arguments follow the work of Candy \cite{C}.
\end{abstract}

\section{Chern-Simons-Dirac equation in one spatial dimension}
We consider the Cauchy problem for the Chern-Simons-Dirac equation in 
one spatial dimension. Let $\psi={}^t(\psi_1(t,x), \psi_2(t,x))
:\re\times \re \to \mathbb{C}^2$ and $\A={}^t(A_0(t,x),A_1(t,x)):
\re\times \re \to \re^2$ be the unknown functions subject to the 
Cauchy problem for the 
Chern-Simons-Dirac system:
\begin{equation} \label{eqn;CSD}
\left\{
\begin{aligned}
    &i\gm^\mu (\pt_\mu-iA_{\mu}) \psi =m\psi,\quad
        &t>0, x\in \re,\\
    &\pt_t A_1-\pt_x A_0=\psi^* \al \psi,\quad  
        &t>0, x\in \re,\\
    &\pt_t A_0-\pt_x A_1=0,\quad   
        &t>0, x\in \re,\\
    &\psi(0,x)=\psi_0(x), \quad
     \A(0,x)=\A_{0}(x),  &x\in \re,
\end{aligned}
\right.
\end{equation}
where 
$m>0$ is a constant, $\mu=0,1$
and $\gm^\mu$ are Dirac matrices
\begin{equation}
\gm^0=
\begin{pmatrix}
1 &0 \\
0 &-1
\end{pmatrix},
\quad\gm^1=
\begin{pmatrix}
0 &1 \\
-1 &0
\end{pmatrix}.
\end{equation}
The initial data $\psi_{0}={}^t(\psi_{10}(x), \psi_{20}(x))$ and
$\A_{0}={}^t(A_{00}(t,x), A_{10}(t,x))$. 
The coupling matrix $\al$ can be chosen either
 $\gm^0, -i\gm^1$ or identity matrix $I$. 
The conjugate ${}^{*}$ means 
$\psi^{*}={}^{t}\overline{\psi}=(\overline{\psi}_{1},\overline{\psi}_{2})$ 
for $\psi={}^{t}(\psi_{1},\psi_{2})$.

The Chern-Simons gauge theory coupled with the Dirac particles is
described by the Chern-Simons-Dirac Lagrangean density:
$$
 \mathcal{L}=\frac{k}{4}\ep_{\mu,\nu\,\r}F_{\mu,\nu}A_{\r}
         +\psi^*\gm^{\mu}(\pt_{\mu}-iA_{\mu})\psi +m\psi^*\psi,
$$
where $\ep_{\mu,\nu,\r}$ denotes the cyclic change anti-symmetric
tensor.  From physical point of view, $1+2$ dimensional model 
is natural and important in the gauge theory.
On the other hand, it is known that the one spatial dimensional 
model for the Maxwell-Klein-Gordon system is integrable 
and the solution is
 expressed
  by the initial data explicitly  if we 
consider the massless case.
As a natural extension, the one spatial dimensional Chern-Simons-Dirac system may also have 
a favorable mathematical structure 
from the mathematical point of view.   
Unfortunately, it is not 
 the integrable system, we may develop some mathematical 
 theory  on the global well-posedness which has a close nature to the 
 massless Maxwell-Klein-Gordon system.

If we neglect the mass term $m\psi$ in the Dirac part, the system 
maintains a scaling invariant under the transform
\eq{\label{eqn;scaling}
\spl{
   &  \psi_{\lam}(t,x)=\lam \psi(\lam t, \lam x),\\
   &   \A_{\lam}(t,x)=\lam \A(\lam t, \lam x)
}
}
for any $\lam>0$.
Therefore it is natural to consider the well-posedness of the system
in some Hilbert scale where the above scaling left invariant.
Toward this idea, Bournaveas-Candy-Machihara \cite{BCM}
considered the time local well-posedness in the space
$C(I;H^{-1/2+})\times C(I;H^{-1/2+})$, where $H^s=H^s(\re)$ 
is the Sobolev space defined by 
$H^s(\re)=\{f\in \mathcal{S}^*;\<\xi\>^s\hat{f}(\xi)\in L^2(\re)\}$,
where $\hat f$ denotes the Fourier transform of $f$.
We simply denote $H^{-1/2+}=H^{-1/2+\ep}(\re)$ for any small $\ep>0$.

On the other hand, for the global well-posedness, the charge conservation 
law for the Dirac part is used as an a priori estimate:
$$
  \|\psi_1(t)\|_2^2+\|\psi_2(t)\|_2^2
  =\|\psi_1(0)\|_2^2+\|\psi_2(0)\|_2^2.
$$
Hence the global well-posedness is restricted in the class of 
$L^2(\re)$, while the scaling invariant critical space stays in 
$H^{-1/2}(\re)$ in the case \eqref{eqn;CSD}.

In this paper, we choose a different approach to consider the 
time global well-posedness of \eqref{eqn;CSD}.
One simple way is to consider the Banach space to construct 
the time local solution. Thanks to the good structure of the 
system, we may derive the time global a priori estimate in 
$L^p$, where $2\le p\le \infty$ and we may show the 
global existence and well-posedness for the solution 
to \eqref{eqn;CSD} for any large data.  Besides, if we restrict 
ourselves to particular nonlinear coupling cases $\al=\gm^0$ 
or $\al=-i\gm^1$, we may reach much further.
Indeed, if we choose $\al=\gm^0$ or 
$\al=-i\gm^1$ then the nonlinear coupling term
 hold {\it the null structure}
and the well-posedness result is improved.

If the coupling term is given by either of the matrices
$\al=\gm^0$ or $\al=-i\gm^1$, then 
the coupling term is less interactive and 
it satisfies the null condition.
In this case we may have the well-posedness 
result up to the critical Banach scale:
 
\begin{thm}[the null  case]\label{thm;main}
Let $1\le  p\le\infty$ and $\al=\gm^0$ or $\al=-i\gm^1$
 in \eqref{eqn;CSD}.
For any $\psi_0\in L^p(\re)$ and
$\A_0\in  L^p(\re)$, there exists 
 a global  weak solution $(\psi, \A) \in C([0,\infty);L^p)
\times C([0,\infty);L^p)$ to \eqref{eqn;CSD} 
such  that
\begin{enumerate}
\item the solution is unique in 
\eqn{
 \spl{
 &\psi\in C([0,\infty);L^p), \quad
 \A\in C([0,\infty);L^p), \\
 &\psi^{*}\alpha\psi\in L^{p}_{loc}([0,\infty);L^{p}), \quad
 \A\psi\in L^{p}_{loc}([0,\infty);L^{p}).
   }
}
\item The map from the data $(\psi_0,\A_0)$ to 
the solution $(\psi, \A)$ is Lipschitz continuous 
from $L^p\times L^p$ 
to $C([0,\infty);L^p)\times 
C([0,\infty);L^p)$.
\end{enumerate}
\end{thm}

We state next that the nonlinear coupling is given by 
$\al=I$:

\begin{thm}[the non-null case]\label{thm;sub}
Let $1\le p\le\infty$ and $\al=I$ in \eqref{eqn;CSD}.
For any $\psi_0\in L^p(\re)\cap L^{\infty}(\re)$ and
$\A_0\in  L^p(\re)$, there exists 
 a unique global  weak solution 
 $(\psi, \A) \in C([0,\infty);L^p\cap L^{\infty})
\times C([0,\infty);L^p)$ to \eqref{eqn;CSD}.
The map from the data $(\psi_0,\A_0)$ to 
the solution $(\psi, \A)$ is Lipschitz continuous 
from $(L^p\cap L^{\infty})\times L^p$ 
to $C([0,\infty);L^p\cap L^{\infty})\times 
C([0,\infty);L^p)$.
\end{thm}

Theorem \ref{thm;sub} is obtained by noticing that the 
solution of the Dirac part of the system maintains a 
favorable a priori bound.  We call this estimate 
as an intrinsic $L^{\infty}$ estimate and this was 
known by several authors \cite{D}, \cite{B1}, 
\cite{C}.
If the coupling matrix is $I$, then the nonlinear 
coupling does not give no good effect  in the sense 
of the regularity, the result so far is optimal in our 
setting.

\vskip4mm\noindent
{\bf Remark.} The interesting case
 in the above theorem 
is the scaling critical setting $p=1$, where the norm of the space 
maintains the scaling \eqref{eqn;scaling} invariant.
This stands for the space is critical for the local and 
global well-posedness for the system \eqref{eqn;CSD}.
We note that neither the global well-posedness nor 
local well-posedness is yet established 
in the critical Sobolev space $H^{-1/2}(\re)\times H^{-1/2}(\re)$.
It is also worth to remark that there is no any inclusion 
relation between the spaces $L^1(\re)$ and $H^{-1/2}(\re)$,
while if  we consider the case $p>1$ then it holds that $L^p(\re)\subset
H^{-1/2+\ep}(\re)$ where $\ep>0$ is depending on $p$.

The main part of our result is to consider the mass concentration  
phenomena.
Since we treat the scaling critical case, it is required 
to show that  the mass
concentration phenomena; namely for any $\del>0$
$$
  \liminf_{t\to T} \left(\int_{B_{\del}(x_0)}|\psi(t,x)|dx 
  +\int_{B_{\del}(x_0)}|\A(t,x)|dx\right) \ge \ep_0
$$ for  some $\ep_0>0$ for the solution $\psi, \A$
does not occur  in $L^1$.
In order to show this, we derive 
the non-concentration lemma which is possible 
in aid of  the intrinsic $L^{\infty}$ estimate 
for the nonlinear part of the solution.
This idea originally goes back to the $L^{\infty}$
estimate due to Delgado \cite{D}.  Candy \cite{C} extended this idea to the 
nonlinear component for the Dirac part and applied for the global well-posedness of the 1 dimensional Dirac equation
with cubic nonlinearity in $L^2(\re)$.
We extend the intrinsic $L^{\infty}$ estimate 
to the Lebesgue spaces $L^p$, where $1\le p<\infty$. For this estimate, we 
essentially use the null structure of the nonlinear 
couplings for both 
Dirac and gauge equations.

We give the notation here. For any Banach space $X$ with respect to 
$x$ variable (space variable), we denote the Bochner space
$L^p(0,T;X)$ by the abbreviation $L^p_TX$ for  $T>0$ 
(see section 3).  Then the local well-posedness of the system 
can be shown in $L^p$ including the critical case $p=1$
under the smallness condition on the data.
The finite time propagation of the solution and the 
intrinsic $L^{\infty}$ estimate shows the existence of the 
large data solution in time globally.

Before closing this section, we state the each roles for the 
sections below. In section 3, we show the time local solution
for small data by using standard iteration argument.
In section 4, we discuss the finite speed of propagation of 
solutions to the transport equation, which is maybe well-known 
among the researchers of wave equation. 
Thanks to this property, we will remove
the smallness assumption for the data to have time local solution.
In section 5, we show that the solutions never concentrate
anywhere, and we can extend those solutions to global time.

%
\setcounter{equation}{0}
\section{Bilinear estimate in $L^{p}$}
For the proof, we reduce the 
system \eqref{eqn;CSD} into the diagonal form for
 the principal part.
We let $\psi_{\pm}=\psi_1\pm\psi_2$ and 
$A_{\pm}=A_0\mp A_1$ to rewrite for the system
\eqref{eqn;CSD} as
\eq{ \label{eqn;csd}
 \left\{
 \spl{
   &\pt_t\psi_{\pm}
      {\pm}\pt_x\psi_{\pm}=iA_{\mp}\psi_{\pm}-im\psi_{\mp},
      &t\in \re, x\in \re,\\
   &\pt_tA_{\pm}{\pm}\pt_xA_{\pm}=\mp P(\psi_{+},\psi_{-}),
      &t\in \re, x\in \re, \\
   &\quad \psi_{\pm}(0,x)=\psi_{10}(x)\pm \psi_{20}(x),
      &x\in \re,\\
   &\quad \A_{\pm}(0,x)=A_{00}(0)\mp A_{110}(0),
        &x\in \re,
 }
 \right.
}
where 
\begin{equation}\label{coupling}
P(\psi_{+},\psi_{-})=
\begin{cases}
\text{Re}(\psi_{+}\bar{\psi}_{-}^*),\quad &\al=\gm^0, \\
\text{Im}(\psi_{+}\bar{\psi}_{-}^*),\quad &\al=-i\gm^1, \\
(|\psi_+|^2+|\psi_-|^2)/2,\quad &\al=I.
\end{cases}
\end{equation}
We give an explicit formula for the solution for the inhomogeneous 
transport equation. 

\begin{lem}\label{lem;transport-eq}
Let $f_0$ and $F=F(t,x)$ be a locally integrable functions. 
Then the solution to the transport equation
\begin{equation}\label{transport-eq1}
 \left\{
 \begin{split}
     &\pt_tu_\pm \pm\pt_xu_\pm =F_\pm(t,x),\\
     &\quad u_\pm(0)=u_{\pm0}
 \end{split}
 \right.
\end{equation}
is given by 
\begin{equation}\label{transport-eq2}
   u_\pm(t,x)=u_{\pm0}(x\mp t)
     +\int_0^tF_\pm\big(s,x\mp(t-s)\big)ds,
\end{equation}
where the signs $\pm$ are the same in upper case and lower case,
 respectively. 
 Besides if $u_{\pm0}\in L^p$ and $F\in L^1(0,T;L^p)$,
then $u_{\pm}\in C([0,T);L^p)$.
\end{lem}
We set the following the characteristic function on the triangle area.
For $x_0\in\re$ and $R>0$,
\begin{equation*}
\chi_{\Omega_R(x_0)}(t,x)=
\begin{cases}
1, &(t,x)\in\Omega_R(x_0)\\
0, &(t,x)\notin\Omega_R(x_0)
\end{cases}
\end{equation*}
where
\begin{equation*}
\Omega_R(x_0)=\{(t,x)\in\re^2:|x-x_0|\le R-t, 0\le t\le R\}.
\end{equation*}
We also set the characteristic function on the interval,
\begin{equation*}
\chi_{I_R(x_0)}(x)=
\begin{cases}
1, &x\in I_R(x_0)\\
0, &x\notin I_R(x_0)
\end{cases}
\end{equation*}
where
\begin{equation*}
I_R(x_0)=\{x\in\re:|x-x_0|\le R\}.
\end{equation*}
From \eqref{transport-eq2} and the influence region, 
we can put the characteristic function into the integral. 
For either cases $\pm$ in \eqref{transport-eq2},
\begin{align}
&\chi_{\Omega_R(x_0)}(t,x)
u_\pm(t,x) \notag\\
&=\chi_{\Omega_R(x_0)}(t,x)u_{\pm0}(x\mp t)
     +\chi_{\Omega_R(x_0)}(t,x)\int_0^tF_\pm\big(s,x\mp(t-s)\big)ds\notag\\
&=\chi_{\Omega_R(x_0)}(t,x)(\chi_{I_R(x_0)}u_{\pm0})(x\mp t)
     +\chi_{\Omega_R(x_0)}(t,x)\int_0^t(\chi_{\Omega_R(x_0)}F_\pm)
\big(s,x\mp(t-s)\big)ds \label{influence1}
\end{align}

We show a bilinear estimate for the inhomogeneous 
transport equation, which plays an important role for 
the contraction mapping method.
The estimate below with $p=2$ was shown in \cite{Ma}.
The other cases $0<p\le\infty$ is shown in the same way.

\vskip3mm
\begin{lem}\label{lem;bilinear}
For $0< p\le \infty$, let $F_{\pm}\in L^p(\re)$, with initial
data $u_{\pm0}\in L^p$. 
Then  the solutions $u_+,u_-$  to the  transport equation
\begin{equation*}
 \left\{
 \begin{split}
     &\pt_tu_\pm \pm\pt_xu_\pm =F_\pm(t,x),\\
     &\quad u_\pm(0)=u_{\pm0}
 \end{split}
 \right.
\end{equation*}
are subject to the following estimate: For $T>0$,
\begin{equation} \label{eqn;bilinear}
\|u_+u_-\|_{L_T^pL^p}
  \le\Big(\frac12\Big)^{\frac1{p}}
  \Big(\|u_{+0}\|_{L^p}+\int_0^T\|F_+(s)\|_{L^p}ds\Big)\ 
  \Big(\|u_{-0}\|_{L^p}+\int_0^T\|F_-(s)\|_{L^p}ds\Big).
\end{equation}
From \eqref{influence1}, we also have the bilinear estimates
on the restricted domains,
 \begin{align}\label{influence2}
\|\chi_{\Omega_R(x_0)}u_+u_-\|_{L_T^pL^p}
\le\Big(\frac12\Big)^{\frac1{p}}
&(\|\chi_{I_R(x_0)}u_{+0}\|_{L^p}+\int_0^T\|(\chi_{\Omega_R(x_0)}F_{+})(s)\|_{L^p}ds)\\
&\cdot(\|\chi_{I_R(x_0)}u_{-0}\|_{L^p}+\int_0^T\|(\chi_{\Omega_R(x_0)}F_{-})(s)\|_{L^p}ds).\notag
\end{align}
\end{lem}

\vskip2mm
\begin{prf}{Lemma \ref{lem;bilinear}}
From Lemma \ref{lem;transport-eq}, we have for $u_+$ and 
$u_-$ that
\begin{equation*}
   u_\pm(t,x)=u_{\pm0}(x\mp t)
     +\int_0^tF(s,x\mp t\pm s)ds.
\end{equation*}
Changing  the variable $x-t=y, x+t=z$, we see
\begin{align*}
  \|u_{+0}(x-t)u_{-0}(x+t)\|_{L^p_TL^p}
    \le &\Big(\frac12\Big)^{\frac1{p}}
       \|u_{+0}\|_{L^p}\|u_{-0}\|_{L^p}.
\end{align*}
We remark that the above estimate holds for 
$p=\infty$ with   
$(\frac12)^{\frac1{\infty}}=1$.
Analogously
\begin{align*}
 \left\|u_{+0}(x-t)
    \int_0^tF_-(s,x+t-s)ds\right\|_{L^p_TL^p}
  &\le\int_0^T\|u_{+0}(x-t)F_-(s,x+t-s)\|_{L^p_TL^p}ds \\
  &=\Big(\frac12\Big)^{\frac1{p}}\|u_{+0}\|_{L^p}
    \int_0^T\|F_-(s)  \|_{L^p}ds.
\end{align*}
The estimate for the term 
$u_{-0}(x+t)\int_0^t F_+(s,x-t+s)ds$
is similarly obtained.  For the product of the 
external force term, we see in a similar way as
\begin{align*}
   &\left\|\int_0^tF_+(s,x-t+s)ds
     \int_0^tF_-(s,x+t-s)ds\right\|_{L^p_TL^p}
\\
   &\le\int_0^T
       \int_0^T
           \|F_+(s,x-t+s)F_-(s',x+t-s')\|_{L^p_TL^p}dsds' \\
   &=\Big(\frac12\Big)^{\frac1{p}}
       \iint_{(0,T)\times (0,T)}
       \left(
       \int_{\re}
       \int_{z_-\le z_+\le z_-+\sqrt2T}
         |F_+(s,s-z_+)F_-(s',z_--s')|^pdz_+dz_- 
       \right)^{1/p}dsds' 
\\
   &\le\Big(\frac12\Big)^{\frac1{p}}
       \iint_{(0,T)\times (0,T)}
       \left(   
       \int_{\re}
         |F_+(s,s-z_+|^pdz_+
       \int_{\re}
         F_-(s',z_--s')|^pdz_- 
       \right)^{1/p}dsds' 
\\
   &=\Big(\frac12\Big)^{\frac1{p}}
          \int_0^T\|F_+(s)\|_{L^p}ds
          \int_0^T\|F_-(s')\|_{L^p}ds'.
\end{align*}
Those estimates imply the desired bound \eqref{eqn;bilinear}. 

We apply these estimates above to \eqref{influence1} to 
derive \eqref{influence2}.
\end{prf}

%
\section{finite speed of propagation}
We introduce the well known fact, that is, the finite speed of propagation
of wave and transport equations.

\begin{lem}\label{lem;L-infty-finite speed}
Suppose $(u_{\pm},v_{\pm})\in L^{\infty}_TL^{\infty}\times L^{\infty}_TL^{\infty}$ 
be a solution to \eqref{transport-eq1}, so satisfying
\begin{align} 
u_{\pm}(t,x)&=u_{\pm0}(x\mp t)
    +\int_0^tF(u_{\pm},v_{\pm})(s,x\mp(t-s))ds,
     \label{eqn;finite-spead-eq1} \\
v_{\pm}(t,x)&=v_{\pm0}(x\mp t)
    +\int_0^tG(u_{\pm},v_{\pm})(s,x\mp(t-s))ds,
    \label{eqn;finite-spead-eq2}
\end{align}
respectively, where $F(u_{\pm},v_{\pm})$ and $G(u_{\pm},v_{\pm})$ are
any quadratic form like $uu,uv,vv$, where $u$ and $v$ stand for any
$u_{+},u_{-}$ and $v_{+},v_{-}$ respectively.
If for some $x_0\in \re$ and  $R>0$,
\begin{equation*}
u_{\pm0}(x)=v_{\pm0}(x)=0, \qquad |x-x_0|<R,
\end{equation*}
then for any $0<t<R$,
\begin{equation*}
u_{\pm}(t,x)=v_{\pm}(t,x)=0, \qquad |x-x_0|<R-t.
\end{equation*}
In particular if $u_{\pm0}\equiv 0$ on 
$|x-x_0|<R$, then $u_{\pm}(t,x)\equiv 0$ in 
$|x-x_0|<R-t$.
\end{lem}

\vskip3mm
\begin{prf}{Lemma \ref{lem;L-infty-finite speed}}
Fix $x_{0}$ and $R$. 
For any $0<t<R$, we define the following norm which is the supremum
on the triangle $\Delta=\Delta(x_{0},R)$,
\begin{equation*}
\|f(t)\|_{\Delta}:=\sup_{\{x:|x-x_0|<R-t\}}|f(t,x)|.
\end{equation*}
For any $0<s<t$, we have
\begin{align*}
&\{|x-t+s-x_0|<R-s\}\supset\{|x-x_0|<R-t\}, \\
&\{|x+t-s-x_0|<R-s\}\supset\{|x-x_0|<R-t\}.
\end{align*}
We consider the integral equation \eqref{eqn;finite-spead-eq1}
with $F(u_{\pm},v_{\pm})=u_{\pm}v_{\pm}$ but
the initial data disappear from the assumption,
\begin{align*}
\|u_{\pm}(t)\|_{\Delta}
&\le \sup_{\{x:|x-x_0|<R-t\}} 
        \int_0^t
          |(u_{\pm}v_{\pm})(s,x-t+s)|ds\\
&\le\int_0^t
     \sup_{\{x:|x-t+s-x_0|<t_0-s\}}|(u_{\pm}v_{\pm})(s,x-t+s)|ds \\
&\le \|v_{\pm}\|_{L^{\infty}_TL^{\infty}}
     \int_0^t\|u_{\pm}(s)\|_{\Delta}ds.
\end{align*}
Then Gronwall's inequality implies $\|u_{\pm}(t)\|_{\Delta}=0$. 
Other cases of the quadratic coupling follow similarly. 
\end{prf}

From this argument of the proof of lemma, we may say if two solutions 
$(u_{\pm}^{(1)},v_{\pm}^{(1)}), (u_{\pm}^{(2)},v_{\pm}^{(2)})$ are 
coincide when initial time  
$u_{\pm0}^{(1)}=u_{\pm0}^{(2)}, \ v_{\pm0}^{(1)}=v_{\pm0}^{(2)}$ 
where $|x-x_{0}|<R$, then the solutions coincide
$u_{\pm}^{(1)}=u_{\pm}^{(2)}, \ v_{\pm}^{(1)}=v_{\pm}^{(2)}$ 
where $|x-x_{0}|<R-t$. 
We check this for the case 
$F(u_{\pm},v_{\pm})=G(u_{\pm},v_{\pm})=u_{\pm}v_{\pm}$ only here.
\begin{align*}
\|&u_{\pm}^{(1)}(t)-u_{\pm}^{(2)}(t)\|_{\Delta}
\le \sup_{\{x:|x-x_0|<R-t\}} 
        \int_0^t
          |(u_{\pm}^{(1)}v_{\pm}^{(1)}-u_{\pm}^{(2)}v_{\pm}^{(2)})(s,x-t+s)|ds\\
&\le \|u_{\pm}^{(1)}\|_{L^{\infty}_TL^{\infty}}
     \int_0^t\|v_{\pm}^{(1)}(s)-v_{\pm}^{(2)}\|_{\Delta}ds
+     \|v_{\pm}^{(2)}\|_{L^{\infty}_TL^{\infty}}
     \int_0^t\|u_{\pm}^{(1)}(s)-u_{\pm}^{(2)}(s)\|_{\Delta}ds.
\end{align*}
The same estimate gives $\|v_{\pm}^{(1)}(t)-v_{\pm}^{(2)}(t)\|_{\Delta}$
is bounded by the same right hand side of this. 
Add up of these two estimates and apply 
Gronwall's inequality to get the conclusion. 


%

\section{Intrinsic $L^{\infty}$ estimate for Dirac part}
\setcounter{equation}{0} 
In this section, we show the $L^{\infty}$ estimate for 
the solution of the Dirac part of the system.
We  consider the diagonal form of the system 
\eqref{eqn;csd}, \eqref{coupling}.

We introduce the decomposition for the 
solution of the Dirac equation $\psi_{\pm}$
such as 
\eqn{
   \psi_{\pm}=\psi_{L\pm}+\psi_{N\pm},   
}
where each component of solution satisfies 
the following equations, respectively:
\eq{ \label{eqn109-2}
   \left\{
 \spl{
   &\pt_t\psi_{L\pm}
      \pm \pt_x\psi_{L\pm}=iA_{\mp}\psi_{L\pm},
       &t\in \re, x\in \re,  \\
   &\pt_t\psi_{N\pm} 
      \pm \pt_x\psi_{N\pm}=iA_{\mp}\psi_{N\pm}
        -im\psi_{\mp},
       &t\in \re, x\in \re, \\
   &\pt_tA_{\pm}{\pm}\pt_xA_{\pm}=-P(\psi_{\pm},\psi_{\mp}),
      &t\in \re, x\in \re, \\
   &\quad \psi_{L\pm}(0)=\psi_{{\pm}0},
    \quad \psi_{N\pm}(0)=0,
     }
    \right.
   }
where the nonlinear coupling is given by any of 
\eqref{coupling}.
We then show that the nonlinear coupling part 
$\psi_{\pm}^N$ satisfies {\it the intrinsic $L^{\infty}$ bound}
as follows:

\vskip3mm
\begin{prop}
[Intrinsic $L^{\infty}$ estimate]\label{prop;intrinsic-bound}
Let $\psi_{L\pm}(t,x)$, $\psi_{N\pm}(t,x)$ 
be smooth functions solving the above system
\eqref{eqn109-2}.  Then it satisfies that
\eq{\label{eqn;good-equal}
   |\psi_{L\pm}(t,x)|=|\psi_{\pm0}(x\mp t)|
}
and 
\eqn{ \label{eqn;L-infty}
     \|\psi_{N+}(t)\|_{L^p\cap L^{\infty}}
        +\|\psi_{N-}(t)\|_{L^p\cap L^{\infty}} 
     \le  m\big(\|\psi_{+0}\|_{L^p}
              +\|\psi_{-0}\|_{L^p}\big)(e^{mt}+t-1).
    }
\end{prop}
\begin{cor}
Let $\psi_{\pm}(t,x)$
be smooth functions solving the system 
\eqref{eqn;csd}, \eqref{coupling}.  Then it satisfies that
\begin{equation}
\|\psi_{\pm}(t)\|_{L^p}
     \le  C\big(\|\psi_{+0}\|_{L^p}
              +\|\psi_{-0}\|_{L^p}\big)(e^{mt}+t).
\end{equation}
\end{cor}

\begin{prf}{Proposition \ref{prop;intrinsic-bound}}
Suppose that $(\psi_{L\pm},\psi_{N\pm})$ is smooth 
function and satisfies the above system: 
Then noticing
\begin{align*}
   \pt_t|\psi_{L+}|^2+\pt_x|\psi_{L+}|^2
    =2\text{Re}\big(\bar{\psi}_{L+}(\pt_t\psi_{L+}+\pt_x\psi_{L+})\big)
    =2\text{Re}\big(iA_-|\psi_{L+}|^2\big)=0,
\end{align*}
we deduce 
\begin{equation}   
   |\psi_{L+}(t,x)|=|\psi_{+0}(x-t)|.
\end{equation}
Similarly
\begin{align*}
    \pt_t|\psi_{N+}|^2+\pt_x|\psi_{N+}|^2
     =2\text{Re}\big(iA_-|\psi_{N+}|^2-im\psi_-
      \bar{\psi}_{N+}\big)
     =2m\text{Im}\big(\psi_- \bar{\psi}_{N+}\big).
\end{align*}
While
\begin{align*}
    \pt_t|\psi_{N+}|^2+\pt_x|\psi_{N+}|^2
      =2|\psi_{N+}|(\pt_t|\psi_{N+}|+\pt_x|\psi_{N+}|).
\end{align*}
Hence  if  $|\psi_{N+}|\neq 0$, then  
\begin{equation}
    \pt_t|\psi_{N+}|+\pt_x|\psi_{N+}|
      =m\text{Im}\left(\frac{\psi_- \bar{\psi}_{N+}}{|\psi_{N+}|}\right).
\end{equation}
Form Lemma \ref{lem;transport-eq}, it follows that
\begin{equation}\label{eqn109-4}
     |\psi_{N+}(t,x)|
     =m\int_0^t\text{Im}\left(\frac{\psi_-\bar{\psi}_{N+}}{|\psi_{N+}|}
     \right)(s,x-t+s)ds.
\end{equation}
Analogously for $\psi_-$,
\begin{align*}
   |\psi_{L-}(t,x)|
    =&\, |\psi_{-0}(x+t)|,  \\
   |\psi_{N-}(t,x)|
    =&\, m\int_0^t\text{Im}
       \left(\frac{\psi_+ \bar{\psi}_{N-}}
                  {|\psi_{N-}|}
       \right)(s,x+t-s)ds.
\end{align*}
By the H\"older inequality, 
\begin{align*}
   |\psi_{N+}(t,x)|
   &\le m\int_0^t|\psi_-(s,x-t+s)|ds \\
   &\le m\int_0^t|\psi_{-0}(2s+x-t)|+|\psi_{N-}(s,x-t+s)|ds \\
   &\le mt^{1/p'}\|\psi_{-0}\|_{L^p}+m\int_0^t\|\psi_{N-}(s)\|_{L^\infty_x}ds,
\end{align*}
where $\frac1{p}+\frac1{p'}=1$. 
Applying similar estimate for $\psi_-$,  we 
derive that
\begin{equation} 
  \begin{split}
    \|\psi_{N+}(t)\|_{L^\infty_x}+\|\psi_{N-}(t)\|_{L^\infty_x}
    \le & mt^{1/p'}(\|\psi_{+0}\|_{L^p}+\|\psi_{-0}\|_{L^p}) \\
        & +m\int_0^t\|\psi_{N+}(s)\|_{L^\infty_x}+\|\psi_{N-}(s)\|_{L^\infty_x}ds.
  \end{split}
\end{equation}
Gronwall's inequality now implies that
\begin{equation}\label{eqn109-5}
     \|\psi_{N+}(t)\|_{L^\infty_x}+\|\psi_{N-}(t)\|_{L^\infty_x}
     \le  m(\|\psi_{+0}\|_{L^p}+\|\psi_{-0}\|_{L^p})(e^{mt}+t-1).
\end{equation}
From \eqref{eqn109-4},  we also have for $1\le p <\infty$ that 
\begin{align*}
    \|\psi_{N+}(t)\|_{L^p}
    &\le
     mt\|\psi_{-0}\|_{L^p}+\int_0^t\|\psi_{N-}(s)\|_{L^p}ds
\end{align*}
and hence we obtain 
\begin{equation}
 \begin{split}
    &\|\psi_{N+}(t)\|_{L^p_x}+\|\psi_{N-}(t)\|_{L^p_x} \\
     &\le m t(\|\psi_{+0}\|_{L^p}+\|\psi_{-0}\|_{L^p})
          +\int_0^t\|\psi_{N+}(s)\|_{L^p_x}+\|\psi_{N-}(s)\|_{L^p_x}ds.
 \end{split}
\end{equation}
We again conclude that
\begin{equation}\label{eqn109-6}
       \|\psi_{N+}(t)\|_{L^p_x}+\|\psi_{N-}(t)\|_{L^p_x}
        \le m(\|\psi_{+0}\|_{L^p}+\|\psi_{-0}\|_{L^p})(e^{t}-1).
\end{equation}
Combining \eqref{eqn109-5} and \eqref{eqn109-6}, we obtain the 
desired estimate.
\end{prf}

We should note that the above estimates in 
Proposition \ref{prop;intrinsic-bound} available with any 
nonlinear coupling with the Chern-Simons gauge part, 
$\alpha=\gamma^{0},-i\gamma^{1}$ and $I$.

%

\section{Time global well-posedness  
for the null case}
\setcounter{equation}{0}
In this section, we give the proof of Theorem \ref{thm;main}.
Here we consider the null 
interaction case $\al=\gm^0$
i.e.,
\eq{ \label{eqn;CSD-null}
 \left\{
 \spl{
   &\pt_t\psi_{\pm}
      \pm\pt_x\psi_{\pm}
          =iA_{\mp}\psi_{\pm}
          -im\psi_{\mp}, &\quad t\in \re, x\in \re,\\
   &\pt_tA_{\pm}\pm\pt_xA_{\pm}
         =\mp\text{Re}(\psi_{\pm}\bar{\psi}_{\mp}),
           &\quad t\in \re, x\in \re,\\
   &(\psi_{\pm},A_{\pm})|_{t=0}  
     =(\psi_{\pm0},A_{\pm0}).
 }
 \right.
}
The other case $\al=-i\gm^1$ is considered as an equivalent form and can be estimated in a similar manner.
We give the proof for the subcritical case $L^{p}, p>1$ and the critical 
case $L^{1}$ separately.

\subsection{Null and subcritical case}
Our aim in this subsection is the following,
\begin{thm}\label{prop1}
Let $\al=\gm^{0}$ or $-i\gm^{1}$. Let $1< p\le \infty$. 
For any $(\psi_{\pm0}, A_{\pm0})\in L^p\times L^p$,  
there exsists a global weak solution 
$(\psi_{\pm},A_{\pm})\in C([0,\infty);L^p)\times C([0,\infty);L^p)$
to \eqref{eqn;CSD-null} such that the solution is unique in 
\eq{\label{unique1}
 \spl{
 &\psi_{\pm}\in C([0,\infty);L^p), \quad
 A_{\pm}\in C([0,\infty);L^p), \\
 & A_{\mp}\psi_{\pm}\in L^{p}_{loc}(0,\infty;L^{p}), \quad
 \psi_{\pm}\bar{\psi}_{\mp}\in L^{p}_{loc}(0,\infty;L^{p}).
   }
}
The map from the data $(\psi_{\pm0},A_{\pm0})$ to 
the solution $(\psi_{\pm}, A_{\pm})$ is Lipschitz continuous 
from $L^p\times L^p$ 
to $C([0,\infty);L^p)\times 
C([0,\infty);L^p)$.
\end{thm}

\begin{prf}{Theorem \ref{prop1}}
We set $1<p\le\infty$. 
We consider the following recurrence of successive approximation of the 
solution.
Let $\{\psi^{(n)}_{\pm},A^{(n)}_{\pm}\}_{n=1,2,\ldots}$ 
solve
\eq{ \label{csd-null-sequence}
 \left\{
 \spl{
   &\pt_t\psi_{\pm}^{(n+1)}
      \pm\pt_x\psi_{\pm}^{(n+1)}=iA_{\mp}^{(n)}\psi_{\pm}^{(n)}
      -im\psi_{\mp}^{(n)},\\
   &\pt_tA_{\pm}^{(n+1)}\pm\pt_xA_{\pm}^{(n+1)}
     =\mp\text{Re}(\psi_{\pm}^{(n)}{\bar{\psi}}_{\mp}^{(n)}),\\
   &(\psi^{(n)}_{\pm},A^{(n)}_{\pm})|_{t=0}  
     =(\psi_{\pm0},A_{\pm0})
 }
 \right.
}
with the first step
$(\psi^{(0)}_{\pm},A^{(0)}_{\pm})=(0,0)$.
Now we consider the integral equation derived from 
Lemma \ref{lem;transport-eq}
\eq{ \label{eqn;int-rec-null}
\left\{
\spl{
  \psi_{\pm}^{(n+1)}(t,x)&=\psi_{\pm0}(x\mp t)
     +i\int_0^t(A_{\mp}^{(n)}\psi_{\pm}^{(n)}
               -m\psi_{\mp}^{(n)})(s,x\mp t\pm s)ds,\\
  A_{\pm}^{(n+1)}(t,x)&=A_{\pm0}(x\mp t)
    \mp \int_0^t(\text{Re}(\psi_{\pm}^{(n)}{\bar{\psi}_\mp}^{(n)}))
       (s,x\mp t\pm s)ds.
  }
  \right.
}
We set
\begin{equation}\label{initial-data}
  M=2(\|\psi_{+0}\|_{L^p}+\|\psi_{-0}\|_{L^p}
      +\|A_{+0}\|_{L^p}+\|A_{-0}\|_{L^p})
\end{equation}
and assume that $M$ is sufficiently small.
We then recall the function space defined by  
\eqn{
  X^p_T
   \equiv  \left.
   \begin{cases}
     &(\phi,\mathbb{B})\in C([0,T);L^p)\times
                 C([0,T);L^p); 
      \qquad 
      \|\phi\|_{L^{\infty}_TL^p}
     +\|\B\|_{L^{\infty}_TL^p}
      \le M
    \end{cases}\right\}.
   }
We show the following bound, 
\begin{equation}\label{each-step-bound11}
  \|\psi_{\pm}^{(n)}\|_{L^{\infty}_TL^p}, \ 
  \|A_{\pm}^{(n)}\|_{L^{\infty}_TL^p}\le M
\end{equation}
for all $n=1,2,\ldots$.
Besides
\begin{equation}\label{each-step-bound12}
  \|\psi_{\pm}^{(n)}A_{\mp}^{(n)}\|_{L^p_TL^p}, \ 
  \|\psi_{\pm}^{(n)}\psi_{\mp}^{(n)}\|_{L^p_TL^p}\le M^2 
\end{equation}
for all $n=1,2,\ldots$. 
We show \eqref{each-step-bound11} and \eqref{each-step-bound12} by induction
with respect to $n$. It is obvious with $n=0$. 
It follows from the first equation of the 
integral form \eqref{eqn;int-rec-null} and H\"older inequality 
that
\begin{align*}
  \|\psi_{+}^{(n+1)}\|_{L^{\infty}_TL^p}
   &\le\|\psi_{+0}\|_{L^p}+\|A_-^{(n)}\psi_+^{(n)}\|_{L^1_TL^p}
    +mT\|\psi_{-}^{(n)}\|_{L^{\infty}_TL^p}\\
   &\le\frac{M}2+T^{1-\frac1p}M^2+mTM\le M.
\end{align*}
The last inequality holds for sufficiently small $T$.
Similarly one may obtain
$\|\psi_{-}^{(n+1)}\|_{L^{\infty}_TL^p}\le M$. 
We also have for small $T$ that  
\begin{align*}
   \|A_{+}^{(n+1)}\|_{L^{\infty}_TL^p}
    \le\|A_{+0}\|_{L^p}+\|\psi_+^{(n)}\psi_-^{(n)}\|_{L^1_TL^p}
    \le\frac{M}2+T^{1-\frac1p}M^2\le M
\end{align*}
and similarly   $\|A_{-}^{(n+1)}\|_{L^{\infty}_TL^p}\le M$. 
According to Lemma \ref{lem;bilinear}, 
we then obtain that
\eqn{
  \spl{
    \|\psi_{+}^{(n+1)}A_-^{(n+1)}\|_{L^{p}_TL^p}
    \le & \Big(\frac12\Big)^{\frac{1}{p}}\Big(\|\psi_{+0}\|_{L^p}
          +\|A_-^{(n)}\psi_+^{(n)}\|_{L^1_TL^p}
          +mT\|\psi_{-}^{(n)}\|_{L^{\infty}_TL^p}\Big) \\
        & \times
         \Big(\|A_{-0}\|_{L^p}+\|\psi_+^{(n)}\psi_-^{(n)}\|_{L^1_TL^p}\Big) \\
    \le& \Big(\frac12\Big)^{\frac{1}{p}}
    \Big(\frac{M}2+T^{1-\frac1p}M^2+mTM\Big)
               \Big(\frac{M}2+T^{1-\frac1p}M^2\Big)
    \le M^{2}
    }
   }
for small $T>0$.
Analogously we obtain 
 $\|\psi_{-}^{(n+1)}A_+^{(n+1)}\|_{L^{p}_TL^p}\le M^{2}$. 
 Furthermore we see that
\begin{align*}
  \|\psi_{+}^{(n+1)}\psi_-^{(n+1)}\|_{L^{p}_TL^p}
  &\le\Big(\frac12\Big)^{\frac{1}{p}}
  (\|\psi_{+0}\|_{L^p}+\|A_-^{(n)}\psi_+^{(n)}\|_{L^1_TL^p}
   +mT\|\psi_{-}^{(n)}\|_{L^{\infty}_TL^p})\\
    &\qquad\qquad\times(\|\psi_{-0}\|_{L^p}+\|A_+^{(n)}\psi_-^{(n)}\|_{L^1_TL^p}
   +mT\|\psi_{+}^{(n)}\|_{L^{\infty}_TL^p})\\
&\le\Big(\frac12\Big)^{\frac{1}{p}}(\frac{M}2+T^{1-\frac{1}{p}}M^2+mTM)^2
\le M^2
\end{align*}
for small $T>0$.
We have obtained \eqref{each-step-bound11} and \eqref{each-step-bound12} for
any $n$.

We next estimate the difference $\psi_{\pm}^{(n+1)}-\psi_{\pm}^{(n)}$
and $A_{\mp}^{(n+1)}-A_{\mp}^{(n)}$.
We set
\begin{align*}
&d((\psi_{\pm},A_{\pm}),(\phi_{\pm},B_{\pm})) \\
=&\|\psi_{\pm}-\phi_{\pm}\|_{L^{\infty}_TL^p}
+\|\psi_{\pm}A_{\mp}-\phi_{\pm}B_{\mp}\|_{L^p_TL^p}
+\|A_{\pm}-B_{\pm}\|_{L^{\infty}_TL^p}
+\|\psi_{\pm}\psi_{\mp}-\phi_{\pm}\phi_{\mp}\|_{L^p_TL^p}.
\end{align*}
We will show 
$d((\psi_{\pm}^{(n+1)},A_{\pm}^{(n+1)}),(\psi_{\pm}^{(n)},A_{\pm}^{(n)}))
\le\frac{1}{2^{n}}$, that is 
\begin{equation}\label{subcritical011}
\begin{split}
&\|\psi_{\pm}^{(n+1)}-\psi_{\pm}^{(n)}\|_{L^{\infty}_TL^p}
+\|\psi_{\pm}^{(n+1)}A_{\mp}^{(n+1)}-\psi_{\pm}^{(n)}A_{\mp}^{(n)}\|
_{L^p_TL^p}\\
&+\|A_{\pm}^{(n+1)}-A_{\pm}^{(n)}\|_{L^{\infty}_TL^p}
+\|\psi_{\pm}^{(n+1)}\psi_{\mp}^{(n+1)}-\psi_{\pm}^{(n)}\psi_{\mp}^{(n)}\|
_{L^p_TL^p}
\le\frac{1}{2^n},
\end{split}
\end{equation}
for small $T>0$ by induction. Indeed
\begin{align*}
\|\psi_{+}^{(n+1)}-\psi_{+}^{(n)}\|_{L^{\infty}_TL^p}
\le  T^{1-\frac{1}{p}}\|A_{-}^{(n)}\psi_{+}^{(n)}-A_-^{(n-1)}\psi_{+}^{(n-1)}\|_{L^p_TL^p}
+mT\|\psi_{-}^{(n)}-\psi_{-}^{(n-1)}\|_{L^{\infty}_TL^p}.
\end{align*}
and
\begin{align*}
\|A_{+}^{(n+1)}-A_{+}^{(n)}\|_{L^{\infty}_TL^p}
\le T^{1-\frac{1}{p}}
\|\psi_{+}^{(n)}\psi_{-}^{(n)}-\psi_+^{(n-1)}\psi_{-}^{(n-1)}\|_{L^p_TL^p}.
\end{align*}
We estimate by the triangle inequality,
\begin{equation*}
\|\psi_{+}^{(n+1)}A_{-}^{(n+1)}-\psi_{+}^{(n)}A_{-}^{(n)}\|
_{L^p_TL^p}
\le
\|\psi_{+}^{(n+1)}(A_{-}^{(n+1)}-A_{-}^{(n)})\|_{L^p_TL^p}
+\|A_{-}^{(n)}(\psi_{+}^{(n+1)}-\psi_{+}^{(n)})\|_{L^p_TL^p}
\end{equation*}
and estimate each term by Lemma \ref{lem;bilinear},
\begin{align*}
&\|\psi_{+}^{(n+1)}(A_{-}^{(n+1)}-A_{-}^{(n)})\|_{L^p_TL^p}\\
&\le\Big(\frac12\Big)^{\frac{1}{p}}(\|\psi_{+0}\|_{L^1}+\|A_-^{(n)}\psi_+^{(n)}\|_{L^p_TL^p}
+mT\|\psi_{-}^{(n)}\|_{L^{\infty}_TL^p})
T^{1-\frac{1}{p}}\|\psi_+^{(n)}\psi_-^{(n)}-\psi_+^{(n-1)}\psi_-^{(n-1)}\|
_{L^p_TL^p}\notag
 \\
&\le\Big(\frac12\Big)^{\frac{1}{p}}(\frac{M}2+M^2+mTM)T^{1-\frac{1}{p}}
\|\psi_+^{(n)}\psi_-^{(n)}-\psi_+^{(n-1)}\psi_-^{(n-1)}\|_{L^p_TL^p},\notag 
\end{align*}
and
\begin{align*}
&\|A_{-}^{(n)}(\psi_{+}^{(n+1)}-\psi_{+}^{(n)})\|_{L^p_TL^p}\\
&\le\Big(\frac12\Big)^{\frac{1}{p}}(\frac{M}2+M^2)(T^{1-\frac{1}{p}}
\|A_{-}^{(n)}\psi_{+}^{(n)}-A_-^{(n-1)}\psi_{+}^{(n-1)}\|_{L^p_TL^p}
+mT\|\psi_{-}^{(n)}-\psi_{-}^{(n-1)}\|_{L^{\infty}_TL^p}). \notag\\
\end{align*}
We estimate by the triangle inequality,
\begin{equation*}
\|\psi_{+}^{(n+1)}\psi_{-}^{(n+1)}-\psi_{+}^{(n)}\psi_{-}^{(n)}\|
_{L^p_TL^p}
\le\|\psi_{+}^{(n+1)}(\psi_{-}^{(n+1)}-\psi_{-}^{(n)})\|_{L^p_TL^p}
+\|\psi_{-}^{(n)}(\psi_{+}^{(n+1)}-\psi_{+}^{(n)})\|_{L^p_TL^p}\end{equation*}
and estimate each term by Lemma \ref{lem;bilinear},
\begin{align*}
&\|\psi_{+}^{(n+1)}(\psi_{-}^{(n+1)}-\psi_{-}^{(n)})\|_{L^p_TL^p} \\
&\le\Big(\frac12\Big)^{\frac{1}{p}}(\frac{M}2+M^2+mTM)(T^{1-\frac{1}{p}}
\|\psi_-^{(n)}A_+^{(n)}-\psi_-^{(n-1)}A_+^{(n-1)})\|_{L^p_TL^p}
+mT\|\psi_+^{(n)}-\psi_+^{(n-1)}\|_{L^{\infty}_TL^p}). \notag\\
&\|\psi_{-}^{(n)}(\psi_{+}^{(n+1)}-\psi_{+}^{(n)})\|_{L^p_TL^p} \\
&\le\Big(\frac12\Big)^{\frac{1}{p}}(\frac{M}2+M^2+mTM)(T^{1-\frac{1}{p}}
\|A_{-}^{(n)}\psi_{+}^{(n)}-A_-^{(n-1)}\psi_{+}^{(n-1)}\|_{L^p_TL^p}
+mT\|\psi_{-}^{(n)}-\psi_{-}^{(n-1)}\|_{L^{\infty}_TL^p}). \notag
\end{align*}
We have obtained \eqref{subcritical011} for any $n$ for small $T$.
Therefore we have obtained the time local existence of solution.
Here we remark the existence time $T$ depends on
the size of norms of initial data 
$\|\psi_{\pm0}\|_{L^{p}}, \|A_{\pm0}\|_{L^{p}}$ only.
From the intrinsic bound, Proposition \ref{prop;intrinsic-bound},
we can extend the solution to the global time.

Next we discuss the uniqueness of solutions. We suppose two solutions
$(\psi_{\pm},A_{\pm}),(\phi_{\pm},B_{\pm})$ with same initial data
to
\eqref{eqn;CSD-null} satisfy the condition \eqref{unique1}. Then there 
exists $M_{1},M_{2}$ such that
\begin{align*}
\|\psi_{\pm}\|_{L^{\infty}_{T}L^{p}}, \ 
\|A_{\pm}\|_{L^{\infty}_{T}L^{p}}\le M_{1}, \quad
\|\psi_{\pm}A_{\mp}\|_{L^{p}_{T}L^{p}}, \ 
\|\psi_{\pm}\psi_{\mp}^{*}\|_{L^{p}_{T}L^{p}}\le M_{1}^{2}, \\
\|\phi_{\pm}\|_{L^{\infty}_{T}L^{p}},
\|B_{\pm}\|_{L^{\infty}_{T}L^{p}}\le M_{2}, \quad
\|\phi_{\pm}B_{\mp}\|_{L^{p}_{T}L^{p}}, \ 
\|\phi_{\pm}\phi_{\mp}^{*}\|_{L^{p}_{T}L^{p}}\le M_{2}^{2}. 
\end{align*}
We then derive the following from the same estimates which have concluded 
\eqref{subcritical011}, 
\begin{equation*}
d((\psi_{\pm},A_{\pm}),(\phi_{\pm},B_{\pm}))\le\frac12
d((\psi_{\pm},A_{\pm}),(\phi_{\pm},B_{\pm}))
\end{equation*}
for small $T>0$, 
which implies the two solutions coincides to each other.

If we set the different initial data 
$(\psi_{\pm0},A_{\pm0})$ and $(\phi_{\pm0},B_{\pm0})$ for 
$(\psi_{\pm},A_{\pm})$ and $(\phi_{\pm},B_{\pm})$ respectively, 
from the same estimates again, 
we have 
\begin{equation*}
d((\psi_{\pm},A_{\pm}),(\phi_{\pm},B_{\pm}))\le\frac12
d((\psi_{\pm},A_{\pm}),(\phi_{\pm},B_{\pm}))
+C(\|\psi_{\pm0}-\phi_{\pm0}\|_{L^{p}}+\|A_{\pm0}-B_{\pm0}\|_{L^{p}})\end{equation*}
for small $T>0$, 
which implies the map from the initial data to the solution is Lipschitz 
continuous.
\end{prf}

\subsection{Null and critical case}\label{null-critical}
Our aim in this subsection is the following,
\begin{thm}\label{thm-null-critical}
Let $\al=\gm^{0}$ or $-i\gm^{1}$. 
For any $(\psi_{\pm0}, A_{\pm0})\in L^1\times L^1$,  
there exsists a global weak solution 
$(\psi_{\pm},A_{\pm})\in C([0,\infty);L^1)\times C([0,\infty);L^1)$
to \eqref{eqn;CSD-null} such that the solution is unique in 
\eq{\label{unique2}
 \spl{
 &\psi_{\pm}\in C([0,\infty);L^1), \quad
 A_{\pm}\in C([0,\infty);L^1), \\
 & A_{\mp}\psi_{\pm}\in L^{1}_{loc}(0,\infty;L^{1}), \quad
 \psi_{\pm}\bar{\psi}_{\mp}\in L^{1}_{loc}(0,\infty;L^{1}).
   }
}
The map from the data $(\psi_{\pm0},A_{\pm0})$ to 
the solution $(\psi_{\pm}, A_{\pm})$ is Lipschitz continuous 
from $L^1\times L^1$ 
to $C([0,\infty);L^1)\times 
C([0,\infty);L^1)$.
\end{thm}

\begin{prf}{Theorem \ref{thm-null-critical}}
A big difference between the following bound and the previous one is 
that the function space is $L^1$ where the solution is scaling 
invariant and there is no room to produce the $T^{\ep}, \ep>0$ 
in this case.
So the proof below is not so straight like the sub-critical case above.
We make some steps for the proof. We show the time local existence of 
solution for small data. We can remove this smallness condition 
of the initial 
data by finite speed of propagation of wave and transport equation. 
We observe the 
non-concentration property of solution which implies the time global 
existence of solution. We show the uniqueness of solution by cut-off 
solution by which we can treat the solution is small.

\subsubsection{Local existence of solution with small size initial data}
In this subsubsection, we will obtain a solution with small data.
We use the same scheme \eqref{csd-null-sequence}, that is 
\eqref{eqn;int-rec-null}, and $M$ of \eqref{initial-data}. 
We consider $M$ is sufficiently small for a while.
We shall show
\begin{equation}\label{each-step-bound1}
\|\psi_{\pm}^{(n)}\|_{L^{\infty}_TL^1}, \ 
\|A_{\pm}^{(n)}\|_{L^{\infty}_TL^1}\le M
\end{equation}
for any $n$ by the induction argument again. We also obtain
the following bounds
\begin{equation}\label{each-step-bound2}
\|\psi_{\pm}^{(n)}A_{\mp}^{(n)}\|_{L^1_TL^1}, \ 
\|\psi_{\pm}^{(n)}\psi_{\mp}^{(n)}\|_{L^1_TL^1}\le M^2 
\end{equation}
for any $n$. We start to estimate
\begin{align*}
\|\psi_{+}^{(n+1)}\|_{L^{\infty}_TL^1}
&\le\|\psi_{+0}\|_{L^1}+\|A_-^{(n)}\psi_+^{(n)}\|_{L^1_TL^1}
+mT\|\psi_{-}^{(n)}\|_{L^{\infty}_TL^1}\\
&\le\frac{M}2+M^2+mTM\le M.
\end{align*}
The last inequality holds for sufficiently small $M$ and $T$.
Similarly we can obtain
$\|\psi_{-}^{(n+1)}\|_{L^{\infty}_TL^1}\le M$. We estimate
\begin{align*}
\|A_{+}^{(n+1)}\|_{L^{\infty}_TL^1}
\le\|A_{+0}\|_{L^1}+\|\psi_+^{(n)}\psi_-^{(n)}\|_{L^1_TL^1}
\le\frac{M}2+M^2\le M.
\end{align*}
Similarly $\|A_{-}^{(n+1)}\|_{L^{\infty}_TL^1}\le M$. 
From Lemma \ref{lem;bilinear}, we have
\begin{align*}
\|\psi_{+}^{(n+1)}A_-^{(n+1)}\|_{L^{1}_TL^1}
&\le\frac12(\|\psi_{+0}\|_{L^1}+\|A_-^{(n)}\psi_+^{(n)}\|_{L^1_TL^1}
+mT\|\psi_{-}^{(n)}\|_{L^{\infty}_TL^1})
(\|A_{-0}\|_{L^1}+\|\psi_+^{(n)}\psi_-^{(n)}\|_{L^1_TL^1}) \\
&\le\frac12(\frac{M}2+M^2+mTM)(\frac{M}2+M^2)\le M^2.
\end{align*}
Similarly $\|\psi_{-}^{(n+1)}A_+^{(n+1)}\|_{L^{1}_TL^1}\le M^{2}$. 
We estimate
\begin{align*}
\|\psi_{+}^{(n+1)}\psi_-^{(n+1)}\|_{L^{1}_TL^1}
&\le\frac12(\|\psi_{+0}\|_{L^1}+\|A_-^{(n)}\psi_+^{(n)}\|_{L^1_TL^1}
+mT\|\psi_{-}^{(n)}\|_{L^{\infty}_TL^1})\\
&\qquad\qquad\times(\|\psi_{-0}\|_{L^1}+\|A_+^{(n)}\psi_-^{(n)}\|_{L^1_TL^1}
+mT\|\psi_{+}^{(n)}\|_{L^{\infty}_TL^1})\\
&\le\frac12(\frac{M}2+M^2+mTM)^2\le M^2.
\end{align*}
We have obtained \eqref{each-step-bound1} and \eqref{each-step-bound2} for
any $n$.

We next estimate the difference $\psi_{\pm}^{(n+1)}-\psi_{\pm}^{(n)}$
and $A_{\mp}^{(n+1)}-A_{\mp}^{(n)}$.
We set
\begin{align}\label{distance}
&d((\psi_{\pm},A_{\pm}),(\phi_{\pm},B_{\pm})) \\
=&\|\psi_{\pm}-\phi_{\pm}\|_{L^{\infty}_TL^p}
+3\|\psi_{\pm}A_{\mp}-\phi_{\pm}B_{\mp}\|_{L^p_TL^p}
+\|A_{\pm}-B_{\pm}\|_{L^{\infty}_TL^p}
+3\|\psi_{\pm}\psi_{\mp}-\phi_{\pm}\phi_{\mp}\|_{L^p_TL^p}. \notag
\end{align}
We will show 
$d((\psi_{\pm}^{(n+1)},A_{\pm}^{(n+1)}),(\psi_{\pm}^{(n)},A_{\pm}^{(n)}))
\le\frac{1}{2^{n}}$, that is 
\begin{equation}\label{each-step-difference}
\begin{split}
&\|\psi_{\pm}^{(n+1)}-\psi_{\pm}^{(n)}\|_{L^{\infty}_TL^1}
+3\|\psi_{\pm}^{(n+1)}A_{\mp}^{(n+1)}-\psi_{\pm}^{(n)}A_{\mp}^{(n)})\|_{L^1_TL^1}\\
&+\|A_{\pm}^{(n+1)}-A_{\pm}^{(n)}\|_{L^{\infty}_TL^1}
+3\|\psi_{\pm}^{(n+1)}\psi_{\mp}^{(n+1)}-\psi_{\pm}^{(n)}\psi_{\mp}^{(n)})\|_{L^1_TL^1}
\le\frac{1}{2^n}
\end{split}
\end{equation}
for any $n$ by the induction argument. Here it comes from a technical reason
to put $3$ before the norms of square terms. 
We show the estimates of upper sign only; $+$ of $\pm$ and $-$ of $\mp$.
\begin{align}\label{difference1}
\|\psi_{+}^{(n+1)}-\psi_{+}^{(n)}\|_{L^{\infty}_TL^1}
\le\|A_{-}^{(n)}\psi_{+}^{(n)}-A_-^{(n-1)}\psi_{+}^{(n-1)}\|_{L^1_TL^1}
+mT\|\psi_{-}^{(n)}-\psi_{-}^{(n-1)}\|_{L^{\infty}_TL^1}.
\end{align}
and
\begin{align}\label{difference2}
\|A_{+}^{(n+1)}-A_{+}^{(n)}\|_{L^{\infty}_TL^1}
\le\|\psi_{+}^{(n)}\psi_{-}^{(n)}-\psi_+^{(n-1)}\psi_{-}^{(n-1)}\|_{L^1_TL^1}.
\end{align}
We use triangle inequality
\begin{align}\label{difference3}
\|\psi_{+}^{(n+1)}A_{-}^{(n+1)}-\psi_{+}^{(n)}A_{-}^{(n)})\|_{L^1_TL^1}
&\le\|\psi_{+}^{(n+1)}(A_{-}^{(n+1)}-A_{-}^{(n)})\|_{L^1_TL^1}
+\|A_-^{(n)}(\psi_{+}^{(n+1)}-\psi_{+}^{(n)})\|_{L^1_TL^1}, \\
\|\psi_{+}^{(n+1)}\psi_{-}^{(n+1)}-\psi_{+}^{(n)}\psi_{-}^{(n)})\|_{L^1_TL^1}
&\le\|\psi_{+}^{(n+1)}(\psi_{-}^{(n+1)}-\psi_{-}^{(n)})\|_{L^1_TL^1}
+\|\psi_{-}^{(n)}(\psi_{+}^{(n+1)}-\psi_{+}^{(n)})\|_{L^1_TL^1}.\label{difference4}
\end{align}
We estimate each term by Lemma \ref{lem;bilinear},
\begin{align}\label{difference5}
&\|\psi_{+}^{(n+1)}(A_{-}^{(n+1)}-A_{-}^{(n)})\|_{L^1_TL^1}\\
&\le\frac12(\|\psi_{+0}\|_{L^1}+\|A_-^{(n)}\psi_+^{(n)}\|_{L^1_TL^1}
+mT\|\psi_{-}^{(n)}\|_{L^{\infty}_TL^1})
\|\psi_+^{(n)}\psi_-^{(n)}-\psi_+^{(n-1)}\psi_-^{(n-1)}\|_{L^1_TL^1}\notag
 \\
&\le\frac12(\frac{M}2+M^2+mTM)
\|\psi_+^{(n)}\psi_-^{(n)}-\psi_+^{(n-1)}\psi_-^{(n-1)}\|_{L^1_TL^1}.\notag 
\end{align}
We make $\frac12(\frac{M}2+M^2+mTM)$ small sufficiently. Similarly
\begin{align}
&\|A_{-}^{(n)}(\psi_{+}^{(n+1)}-\psi_{+}^{(n)})\|_{L^1_TL^1}\label{difference6}\\
&\le\frac12(\frac{M}2+M^2)
(\|A_{-}^{(n)}\psi_{+}^{(n)}-A_-^{(n-1)}\psi_{+}^{(n-1)}\|_{L^1_TL^1}
+mT\|\psi_{-}^{(n)}-\psi_{-}^{(n-1)}\|_{L^{\infty}_TL^1}). \notag\\
&\|\psi_{+}^{(n+1)}(\psi_{-}^{(n+1)}-\psi_{-}^{(n)})\|_{L^1_TL^1}\label{difference7}\\
&\le\frac12(\frac{M}2+M^2+mTM)
(\|\psi_-^{(n)}A_+^{(n)}-\psi_-^{(n-1)}A_+^{(n-1)})\|_{L^1_TL^1}
+mT\|\psi_+^{(n)}-\psi_+^{(n-1)}\|_{L^{\infty}_TL^1}). \notag\\
&\|\psi_{-}^{(n)}(\psi_{+}^{(n+1)}-\psi_{+}^{(n)})\|_{L^1_TL^1}\label{difference8}\\
&\le\frac12(\frac{M}2+M^2+mTM)
(\|A_{-}^{(n)}\psi_{+}^{(n)}-A_-^{(n-1)}\psi_{+}^{(n-1)}\|_{L^1_TL^1}
+mT\|\psi_{-}^{(n)}-\psi_{-}^{(n-1)}\|_{L^{\infty}_TL^1}). \notag
\end{align}
We have obtained \eqref{each-step-difference} for any $n$.
So we have obtained the time local solution with the small data.

We remark if the initial data is small,
we have already obtained that the solution map from initial data
is Lipshitz continuous as well,
\begin{equation*}
d((\psi_{\pm},A_{\pm}),(\phi_{\pm},B_{\pm}))\le C
(\|\psi_{\pm0}-\phi_{\pm0}\|_{L^{1}}+\|A_{\pm0}-B_{\pm0}\|_{L^{1}}).
\end{equation*}

\subsubsection{Local existence of solution with general size initial data}
We remove the smallness condition for initial data by using the
finite speed of propagation of solution.

For any initial data $\psi_{\pm0},A_{\pm0}\in L^{1}$ and any $M>0$, 
there exists $r>0$ such that
\begin{equation}\label{no-shrink}
\sup_{x\in\re}\int_{|x-y|<r}(|\psi_{\pm0}(y)|+|A_{\pm0}(y)|)dy<M.
\end{equation}
With this $r$, we split the real line for the variable $x$ in 
the following two ways,
\begin{equation*}
\re=\bigcup_{j\in\mathbb{Z}}I_{1}^{j}=\bigcup_{j\in\mathbb{Z}}I_{2}^{j}, 
\quad I_{1}^{j}=[jr,(j+1)r], \ I_{2}^{j}=[(j+\frac12)r,(j+\frac32)r]
\end{equation*}
For each $k=1,2, \ j\in\mathbb{Z}$,
if we replace the initial data $\psi_{\pm0}$ and $A_{\pm0}$ by
$\chi_{I_{k}^{j}}\psi_{\pm0}$ and $\chi_{I_{k}^{j}}A_{\pm0}$ respectively,
the data is small and we could find the corresponding solution
$\psi_{\pm}^{k,j}, A_{\pm}^{k,j}$ with the uniform existence time $T$ with
respect to $k$ and $j$ from the argument in the previous 
subsubsection. We consider $(t,x)$ plane in $\re^{2}$ and in which 
the two squares $S^{1,j}$ and $S^{2,j}$,
\begin{align*}
S^{1,j}&=\{(t,x)\in\re^{2}:0\le t<\frac{r}{4}, j+\frac14\le x<j+\frac34\}, \\
S^{2,j}&=\{(t,x)\in\re^{2}:0\le t<\frac{r}{4}, j+\frac34\le x<j+\frac54\}.
\end{align*}
These are all disjointed and satisfy
\begin{equation*}
\bigcup_{k=1,2}\bigcup_{j\in\mathbb{Z}}S^{k,j}
=\{(t,x)\in\re^{2}:0\le t<\frac{r}{4}\}.
\end{equation*}
We also remark these squares are included in the triangles
\begin{align*} 
S^{1,j}&\subset\{(t,x)\in\re^{2}:|x-(j+\frac12)r)|<\frac{r}{2}-t\}
=:\Omega^{1,j}, \\
S^{2,j}&\subset\{(t,x)\in\re^{2}:|x-(j+1)r)|<\frac{r}{2}-t\}
=:\Omega^{2,j}.
\end{align*} 
We define the functions, for $0<t<\min\{T,\frac{r}{2}\}$,
\begin{equation}\label{solution}
\psi_{\pm}(t,x)=\sum_{k=1,2}\sum_{j\in\mathbb{Z}}
\chi_{S^{k,j}}(t,x)\psi_{\pm}^{k,j}(t,x),
\qquad
A_{\pm}(t,x)=\sum_{k=1,2}\sum_{j\in\mathbb{Z}}
\chi_{S^{k,j}}(t,x)A_{\pm}^{k,j}(t,x).
\end{equation}
In the following, we shall show that \eqref{solution} is a solution
of \eqref{eqn;CSD-null}. In order to do this, it is sufficient to prove
that the solution which is restricted in $S^{k,j}$, that is,
\begin{equation*}
\chi_{S^{k,j}}(t,x)\psi_{\pm}(t,x)
=\chi_{S^{k,j}}(t,x)\psi_{\pm}^{k,j}(t,x), \qquad
\chi_{S^{k,j}}(t,x)A_{\pm}(t,x)
=\chi_{S^{k,j}}(t,x)A_{\pm}^{k,j}(t,x)
\end{equation*}
does not influenced
by changing the initial data on the out of the interval $I_{k}^{j}$.
We apply the equality \eqref{influence1} to \eqref{eqn;int-rec-null},
\begin{equation}\label{cut-problem}
\spl{
  &\chi_{\Omega^{k,j}}\psi_{\pm}(t,x) \\
  &=\chi_{\Omega^{k,j}}(\chi_{I_{k}^{j}}\psi_{\pm0})(x\mp t)
     +i\chi_{\Omega^{k,j}}\int_0^t((\chi_{\Omega^{k,j}}A_{\mp})
     (\chi_{\Omega^{k,j}}\psi_{\pm})
               -m(\chi_{\Omega^{k,j}}\psi_{\mp}))(s,x\mp t\pm s)ds,\\
  &\chi_{\Omega^{k,j}}A_{\pm}(t,x)
  =\chi_{\Omega^{k,j}}(\chi_{I^{j}_{k}}A_{\pm0})(x\mp t)
    \mp \chi_{\Omega^{k,j}}\int_0^t(\text{Re}((\chi_{\Omega^{k,j}}\psi_{\pm})
    (\chi_{\Omega^{k,j}}{\psi_\mp^*})))
       (s,x\mp t\pm s)ds.
  }
\end{equation}
If $(t,x)\in S^{k,j}\subset \Omega^{k,j}$, then all 
$\chi_{\Omega^{k,j}}(t,x)=1$ above. We can observe that the solution 
$(\chi_{\Omega^{k,j}}\psi_{\pm},\chi_{\Omega^{k,j}}A_{\pm})=(\psi_{\pm},A_{\pm})$
of \eqref{cut-problem} won't change when the out side of 
$I_{k}^{j}$ for the initial data $(\psi_{\pm0},A_{\pm0})$ changes.

To this end, we simply check that the solution \eqref{solution} 
is in $L^{\infty}_TL^1$.  We estimate 
\begin{align*}
  \|\psi_{\pm}\|_{L^{\infty}_TL^1}
    +\|A_{\pm}\|_{L^{\infty}_TL^1}
   \le &C\sum_{k=1,2}\sum_{j\in\mathbb{Z}}
      (\|\chi_{I_{k}^{j}}\psi_{\pm0}\|_{L^{\infty}_TL^1}
      +\|\chi_{I_{k}^{j}}A_{\pm0}\|_{L^{\infty}_TL^1}) \\
  \le &C\big( 
      \|\psi_{\pm0}\|_{L^1}+\|A_{\pm0}\|_{L^1}\big).
\end{align*}

\subsubsection{Uniqueness of solution}
We now show the uniqueness of the solutions.
We consider two solutions $(\psi_{\pm},A_{\pm})$ and $(\phi_{\pm},B_{\pm})$ 
to \eqref{eqn;CSD-null} satisfying
\begin{equation}\label{unique2}
\psi_{\pm},\ A_{\pm},\ \phi_{\pm},\ B_{\pm}\in L^\infty_TL^1, \qquad
\psi_{\pm}A_{\mp},\  \phi_{\pm}B_{\mp},\ 
\psi_{\pm}\overline{\psi}_{\mp},\ \phi_{\pm}\overline{\phi}_{\mp}
\in L^1_TL^1.
\end{equation}
Fix $x_0\in\re$. For any $M>0$, there is $R$ such that
\begin{equation}\label{unique3}
\|\chi_{\Omega_{R}(x_0)}\psi_{\pm}A_{\mp}\|_{L^1_TL^1}, \ 
\|\chi_{\Omega_{R}(x_0)}\psi_{\pm}\psi_{\mp}\|_{L^1_TL^1},\ 
\|\chi_{\Omega_{R}(x_0)}\phi_{\pm}B_{\mp}\|_{L^1_TL^1}, \ 
\|\chi_{\Omega_{R}(x_0)}\phi_{\pm}\phi_{\mp}\|_{L^1_TL^1}\le M^2.
\end{equation}
We have from \eqref{influence1}
\begin{align*}
\|\chi_{\Omega_{R}(x_0)}\psi_{\pm}\|_{L^{\infty}_TL^1}
\le\|\chi_{I_{R}(x_0)}\psi_{\pm0}\|_{L^1}
+\|\chi_{\Omega_{R}(x_0)}\psi_{\pm}A_{\mp}\|_{L^1_TL^1}
+mT\|\chi_{\Omega_{R}(x_0)}\psi_{\mp}\|_{L^{\infty}_TL^1}.
\end{align*}
We remark the each terms here is finite since \eqref{unique2}.
For the sufficiently small $T$, we have from \eqref{unique3}
\begin{equation}
\|\chi_{\Omega_{R}(x_0)}\psi_{+}\|_{L^{\infty}_TL^1}
+\|\chi_{\Omega_{R}(x_0)}\psi_{-}\|_{L^{\infty}_TL^1}
\le \frac{M}2+4M^2\le M
\end{equation}
by resizing $M$ sufficiently small.
The same estimates for $\chi_{\Omega_{R}(x_0)}\phi_{\pm}$ hold.
We also estimate
\begin{align*}
\|\chi_{\Omega_{R}(x_0)}A_{\pm}\|_{L^{\infty}_TL^1}
\le\|\chi_{I_{R}(x_0)}A_{\pm0}\|_{L^1}
+\|\chi_{\Omega_{R}(x_0)}\psi_{\pm}\psi_{\mp}\|_{L^1_TL^1}
\le\frac{M}2+M^2\le M.
\end{align*}
The same estimates for $\chi_{\Omega_{R}(x_0)}B_{\pm}$ hold.
In the long run, we may think they are small data and small solutions 
to the problem.
From the similar estimates \eqref{difference1}, \eqref{difference2}, 
\eqref{difference3}, \eqref{difference4}, \eqref{difference5}, 
\eqref{difference6}, \eqref{difference7} and  \eqref{difference8}
by using \eqref{influence2},
we have 
\begin{align*}
&d((\chi_{\Omega_{R}(x_0)}\psi_{\pm},\chi_{\Omega_{R}(x_0)}A_{\pm}),
(\chi_{\Omega_{R}(x_0)}\phi_{\pm},\chi_{\Omega_{R}(x_0)}B_{\pm})) \\
&\le\frac12d((\chi_{\Omega_{R}(x_0)}\psi_{\pm},\chi_{\Omega_{R}(x_0)}A_{\pm}),
(\chi_{\Omega_{R}(x_0)}\phi_{\pm},\chi_{\Omega_{R}(x_0)}B_{\pm}))
\end{align*}
where $d((\psi_{\pm},A_{\pm}),(\phi_{\pm},B_{\pm}))$ is defined in
\eqref{distance}.
Therefore we obtained the uniqueness
in the region $\Omega_{R}(x_0)$ at least, but this is sufficient.

\subsubsection{Global existence of solution}
By the previous subsubsection, we have obtained the time local
solution to \eqref{eqn;CSD-null} 
on $[0,T]$. We repeat this argument to have a solution started 
with initial
data $(\psi_{\pm}(T),A_{\pm}(T))$. If we have the condition of
\eqref{no-shrink} for each step, we can derive the global solution.
It is sufficient for this to show the following a priori estimate.
For 
any $T>0$ and any $\ep>0$, there exists 
  $r=r(T,\|\psi_{\pm0}\|_{L^1}
,\|A_{\pm0}\|_{L^1})>0$ such that
\begin{equation}\label{no-shrink2}
\sup_{0<t<T}\sup_{x\in\re}
\int_{|x-y|<r}|\psi_\pm(t,y)|+|A_\pm(t,y)|dy<\varepsilon.
\end{equation}
where $(\psi_{\pm},A_{\pm})$ is the solution for 
\eqref{eqn;CSD-null} in  $C([0,T];L^1)\times 
C([0,T];L^1)$.

We start to prove \eqref{no-shrink2}. 
We use the decomposition $\psi_\pm=\psi_\pm^L+\psi_\pm^N$. From the intrinsic $L^{\infty}$ estimate 
Proposition \ref{prop;intrinsic-bound}, we see that
\begin{equation}
 \label{eqn;intrinsicL-infty}
 \spl{
    \|\psi_{N+}(t)\|_{L^\infty}+\|\psi_{N-}(t)\|_{L^\infty} 
    & \le  m(\|\psi_{+0}\|_{L^1}+\|\psi_{-0}\|_{L^1})(e^{mt}+t-1)\\
    & \le  e^{mT}+T-1,
     }
\end{equation}
Hence by choosing $r>0$ properly small depending on
the right hand side of \eqref{eqn;intrinsicL-infty},
 we have
\begin{equation*}
\sup_{x\in\re}
\int_{|x-y|<r}|\psi_\pm(t,y)|dy
\leq \sup_{x\in\re}
\int_{|x-y|<r}|\psi_{\pm0}(t\mp y)|dy+r(e^{mT}+T-1)<\varepsilon.
\end{equation*}
This provides the non-concentration estimate 
for Dirac part $\psi_{\pm}$. 
We decompose also $A_\pm$ such like
$A_\pm=A_\pm^L+A_\pm^N$
\begin{equation*}
A_\pm^L(t,x)=A_{\pm0}(x\mp t),\qquad
A_\pm^N(t,x)
=\mp\int_0^t\text{Re}(\psi_+\overline{\psi_-})(s,x\mp t\pm s)ds.
\end{equation*}
We need to estimate the integrand of $A^N_+$, that is 
 $\psi_+\psi_-=\psi_{L+}\psi_{L-}+\psi_{L+}\psi_{N-}
+\psi_{N+}\psi_{L-}+\psi_{N+}\psi_{N-}$. 
\eqn{
  \spl{
    \int_0^t|\psi_{L+}\psi_{L-}|(s,x-t+s)ds
      &=|\psi_{+0}(x-t)|\int_0^t|\psi_{-0}(x-t+2s)|ds \\
      &\leq |\psi_{+0}(x-t)|\|\psi_{-0}\|_{L^1}, \\
    \int_0^t|\psi_{L+}\psi_{N-}|(s,x-t+s)ds
      &\le|\psi_{+0}(x-t)|\|\psi_{-}^N\|_{L^\infty}t \\
      &\leq |\psi_{+0}(x-t)|(e^{mT}+T-1)T, \\
    \int_0^t|\psi_{N+}\psi_{L-}|(s,x-t+s)ds
      &\le\|\psi_{+}^N\|_{L^\infty}\int_0^t|\psi_{-0}(x-t+2s)|ds \\
      &\leq (e^{mT}+T-1)\|\psi_{-0}\|_{L^1}, \\
    \int_0^t|\psi_{N+}\psi_{N-}|(s,x-t+s)ds
      &\le\|\psi_{+}^N\|_{L^\infty}\|\psi_{-}^N\|_{L^\infty}T
       \le(e^{mT}+T-1)^2T.
  }
}
Gathering all the estimates, we obtain
\eqn{
  \int_0^t|\psi_+\psi_-|(s,x-t+s)ds
  \leq C|\psi_{+0}(x-t)|+C(T).
  }
Therefore we conclude that 
\begin{align*}
&\sup_{x\in\re}
\int_{|x-y|<r}|A_\pm(t,y)|dy \\
&\leq \sup_{x\in\re}
\int_{|x-y|<r}|A_{\pm0}(t\mp y)|+C|\psi_{+0}(t-y)|
+C|\psi_{-0}(t+y)|+C(T)dy<\varepsilon
\end{align*}
and this shows the estimate for the term of  $A_\pm$.
\end{prf}

%
\setcounter{equation}{0}
\section{Time global well-posedness
 for the non-null case}
In this section we consider the non-null case and
 establish the time local well-posedness for \eqref{eqn;CSD}
 when $\al=I$.
We consider here that 
\eq{ \label{eqn;CSD-non-null}
 \left\{
 \spl{
   &\pt_t\psi_{\pm}
      \pm\pt_x\psi_{\pm}
          =iA_{\mp}\psi_{\pm}
          -im\psi_{\mp}, &\quad t\in \re, x\in \re,\\
   &\pt_tA_{\pm}\pm\pt_xA_{\pm}
         =\frac12 (|\psi_{\pm}|^2+|{\psi}_{\mp}|^2),
           &\quad t\in \re, x\in \re,\\
   &(\psi_{\pm},A_{\pm})|_{t=0}  
     =(\psi_{\pm0},A_{\pm0}).
 }
 \right.
}
\subsection{Non-null and subcritical case}
Our aim in this subsection is the following,

\vskip3mm
\begin{thm}\label{prop3-1}
Let $\al=I$. 
Let $1< p\le \infty$. For any 
$(\psi_{\pm0},A_{\pm0}) \in (L^p\cap L^{\infty})\times L^p$, 
there exists a unique global weak solution 
$(\psi_{\pm},A_\pm)\in C([0,\infty);L^p\cap L^{\infty})\times C([0,\infty);L^p)$ 
to \eqref{eqn;CSD-non-null}. 
The map from the initial data to the solution 
is the Lipschitz continuous from $(L^p\cap L^{\infty})
\times L^p \to C([0,T);L^p\cap L^{\infty})
\times C([0,T);L^p)$.
\end{thm}

\begin{prf}{Theorem \ref{prop3-1}}
We consider the following recurrence of  successive approximation:
Let $\{\psi^{(n)}_{\pm},A^{(n)}_{\pm}\}_{n=1,2,\ldots}$ 
solve
\eq{ \label{eqn;CSD--non-null-rec}
 \left\{
 \spl{
   &\pt_t\psi_{\pm}^{(n+1)}
      \pm\pt_x\psi_{\pm}^{(n+1)}
          =iA_{\mp}^{(n)}\psi_{\pm}^{(n+1)}
      -im\psi_{\mp}^{(n+1)},\\
   &\pt_tA_{\pm}^{(n+1)}\pm\pt_xA_{\pm}^{(n+1)}
         =\frac12 (|\psi_{\pm}^{(n)}|^2+|{\psi}_{\mp}^{(n)}|^2),\\
   &(\psi^{(n)}_{\pm},A^{(n)}_{\pm})|_{t=0}  
     =(\psi_{\pm0},A_{\pm0})
 }
 \right.
}
with the first step
$(\psi^{(1)}_{\pm},A^{(1)}_{\pm})=(\psi_{\pm0},A_{\pm0})$.
We should note that this scheme has different form from 
the former one \eqref{csd-null-sequence}
 not only by the nonlinear coupling 
but also the recurrence suffix on the equation of $\psi$.
The aim for this is to use the calculation for the intrinsic estimates. 
Now then we introduce the integral equation derived from Lemma \ref{lem;transport-eq}:
\begin{equation} \label{eqn;recurrence-scheme-null}
 \left\{ 
 \begin{aligned}
    \psi_{\pm}^{(n+1)}(t,x)&=\psi_{\pm0}(x\mp t)
       +i\int_0^t(A_{\mp}^{(n)}\psi_{\pm}^{(n+1)}
               -m\psi_{\mp}^{(n+1)})(s,x\mp t\pm s)ds,\\
    A_{\pm}^{(n+1)}(t,x)&=A_{\pm0}(x\mp t)
       \mp\frac12  \int_0^t\big(|\psi_{\pm}^{(n)}|^2
                            +|\psi_{\mp}^{(n)}|^2\big)
                            (s,x\mp t\pm s)ds.
 \end{aligned}
 \right.
\end{equation}
Letting
\begin{equation*}
  M=2(\|\psi_{\pm0}\|_{L^p}+\|\psi_{\pm0}\|_{L^{\infty}}
    +\|A_{\pm0}\|_{L^p}).
\end{equation*}
we define  
\eqn{
  X^p_T
   \equiv  \left.
   \begin{cases}
     &(\phi,\mathbb{B})\in C([0,T);L^p)\times
                 C([0,T);L^p); 
      \qquad 
      \|\phi\|_{L^{\infty}L^p}
     +\|\B\|_{L^{\infty}L^p}
      \le M
    \end{cases}\right\}.
   }
It is straightforward to show that $X_T^p$ is a complete 
metric space by the metric induced by the norm:
$$
  \((\phi,\B)\)_{X^p_T}
  \equiv \sup_{t\in[0,T)}
            \|\phi(t)\|_{L^p}
        +\sup_{t\in[0,T)}
           \|\B(t)\|_{L^p}.
$$
For that purpose we show the following estimate for 
the nonlinear coupling:
\begin{lem}\label{lem;nonlinear-bound}
For the sequence 
$\{\psi_{\pm}^{(k)},A_{\pm}^{(k)}\}_{k=1}^{n-1}
  \subset X^p_T$
with $n=2,3,\cdots$, 
it holds that, for small $T>0$
\begin{align}
   & \|\psi_{\pm}^{(n)}\|_{L^{\infty}_T(L^p\cap L^\infty)}\le M,
    \label{eqn;sub-bound-1}\\
   &\||\psi_{\pm}^{(n)}|^2\|_{L^{p}_{T}L^p}
    \le M^2,
    \label{eqn;sub-bound-2}\\
   &\| A_{\pm}^{(n)}\|_{L^{\infty}_TL^p}
     \le M,   
      \label{eqn;sub-bound-3} \\
  &\| A_{\mp}^{(n)}\psi_{\pm}^{(n+1)}\|_{L^{\infty}_TL^p}\le M^2. 
  \label{eqn;sub-bound-4}
\end{align}
\end{lem}
\vskip2mm
\begin{prf}{Lemma \ref{lem;nonlinear-bound}}
Assuming  $\{\psi_{\pm}^{(k)},A_{\pm}^{(k)}\}_{k=1}^{n-1}
\subset X^p_T$ and satisfying the estimate \eqref{eqn;sub-bound-2}
up to $k\le n-1$, 
we decompose $\psi_{\pm}^{(n)}$ into two components
$\psi_{L\pm}^{(n)}+\psi_{N\pm}^{(n)}$ and
it follows from the intrinsic $L^{\infty}$ estimate
(Proposition \ref{prop;intrinsic-bound}) that 
\eq{\label{eqn;psi-N}
   \|\psi_{N\pm}^{(n)}(t)\|_{L^p\cap L^{\infty}}
     \le C e^T\|\psi_{\pm0}\|_{L^p}.
  }
for $0<t<T$. We should note here that the right hand side of 
the recurrence scheme  \eqref{eqn;CSD--non-null-rec}
is linear in $\psi_{\pm}^{(n)}$ which 
is crucial for obtaining the estimate \eqref{eqn;psi-N}.
While we have from
\eqref{eqn;good-equal} that
\eq{ \label{eqn;psi-L}
  \|\psi_{L\pm}^{(n)}\|_{L^{\infty}_T(L^p\cap L^{\infty})}
    =\|\psi_{\pm0}\|_{L^p\cap L^{\infty}}\le \frac14 M.
}
Combining \eqref{eqn;psi-N} and \eqref{eqn;psi-L}
we obtain \eqref{eqn;sub-bound-1},
and which implies \eqref{eqn;sub-bound-2} by H\"older inequality. 
Since $A_{\pm}^{(n)}$ is given by the
solution formula,  we compute from \eqref{eqn;sub-bound-1} that
\eqn{
  \spl{
    \| A_{\pm}^{(n+1)}\|_{L^{\infty}_TL^p}
      \le & \|A_{\pm0}\|_{L^p}
         +\frac12 \int_0^t 
             \big(\||\psi_{+}|^2\|_{L^{\infty}L^p}
                 +\||\psi_{-}|^2\|_{L^{\infty}L^p}\big)ds\\
      \le & \frac12M + CTM^2
      \le M
  }
}
for small $T>0$. \eqref{eqn;sub-bound-1} and \eqref{eqn;sub-bound-3} 
imply \eqref{eqn;sub-bound-4}
\end{prf}

We now show that the recurrence scheme converges to the 
solution of the Chern-Simons-Dirac system:

\vskip3mm
\begin{prop}\label{prop;convergence} Let $\{\psi_{\pm}^{n}\}_{n=1}^{\infty}$
and $\{A_{\pm}^{(n)}\}_{n=1}^{\infty}$ 
defined by recurrence scheme \eqref{eqn;CSD--non-null-rec}.  Then there exists $\psi_{\pm}$,
$A_{\pm}\in C([0,\infty);L^p)$ such that 
\eqn{
\spl{
  &\psi_{\pm}^{(n)} \to \psi_{\pm}, \qquad n\to \infty,\\
  &A_{\pm}^{(n)}\to A_{\pm}, \qquad n\to \infty
}
}
in $C([0,T); L^p)$ and $(\psi_{\pm}, A_{\pm})$ solves the system 
\eqref{eqn;CSD-non-null}.
\end{prop}

\begin{prf}{Proposition \ref{prop;convergence}}
We derive the estimates for the difference
 $\psi_{\pm}^{(n+1)}-\psi_{\pm}^{(n)}$
and $A_{\mp}^{(n+1)}-A_{\mp}^{(n)}$ 
as 
\begin{equation}\label{subcritical01}
\begin{split}
   &\|\psi_{\pm}^{(n+1)}-\psi_{\pm}^{(n)}\|_{L^{\infty}_TL^p}
    +\|A_{\pm}^{(n+1)}-A_{\pm}^{(n)}\|_{L^{\infty}_TL^p}
   \le\frac{1}{2^n}.
\end{split}
\end{equation}
and the sequence is concluded as the Cauchy sequence.
Taking $L^p$ norm to 
\eqn{
\spl{
 \psi_{\pm}^{(n+1)}-\psi_{\pm}^{(n)}
  = &  \int_0^t \big(iA_{\mp}^{(n)}\psi_{\pm}^{(n+1)}
                  -iA_{\mp}^{(n-1)}\psi_{\pm}^{(n)}\big)ds \\
    & +\int_0^t \big(im\psi_{\mp}^{(n+1)}
                  -im\psi_{\mp}^{(n)}\big)ds, 
    }
   }
we see by H\"older inequality that 
\begin{equation}\label{eqn;difference1}
  \begin{split}
   \|\psi_{+}^{(n+1)}-\psi_{+}^{(n)}\|_{L^{\infty}_TL^p}
    \le & T^{1-\frac{1}{p}}
         \|A_{-}^{(n)}\psi_{+}^{(n+1)}-A_-^{(n-1)}\psi_{+}^{(n)}\|_{L^p_TL^p}
        +mT\|\psi_{-}^{(n+1)}-\psi_{-}^{(n)}\|_{L^{\infty}_TL^p}\\
    \le &  T^{1-\frac{1}{p}}(
         \|A_{-}^{(n)}\big(
   \psi_{+}^{(n+1)}-\psi_{+}^{(n)}\big)\|_{L^p_TL^p}
+\|\big(A_{-}^{(n)}-A_-^{(n-1)}\big) \psi_{+}^{(n)}\|_{L^p_TL^p}) \\
& +mT\|\psi_{-}^{(n+1)}-\psi_{-}^{(n)}\|_{L^{\infty}_TL^p}
  \end{split}
\end{equation}
for small $T>0$. 
The first term of the right hand side of \eqref{eqn;difference1} is estimated as 
follows. By Lemma \ref{lem;bilinear} and \eqref{eqn;sub-bound-2},
\eq{\label{eqn000}
\spl{
  &\|A_{-}^{(n)}(\psi_{+}^{(n+1)}-\psi_{+}^{(n)})\|_{L^p_TL^p}\\
   \le& \Big(\frac12\Big)^{\frac{1}{p}}\left(\frac{1}2\|A_0\|_p+M^2\right)
       T^{1-\frac{1}{p}}\left(
        \|A_{-}^{(n)}\psi_{+}^{(n+1)}-A_-^{(n-1)}\psi_{+}^{(n)}\|_{L^p_TL^p}
      +m\|\psi_{-}^{(n+1)}-\psi_{-}^{(n)}\|_{L^{p}_TL^p}\right)
}}
and then separating 
\begin{equation}\label{eqn001}
 \|A_{-}^{(n)}\psi_{+}^{(n+1)}-A_-^{(n-1)}\psi_{+}^{(n)}\|_{L^p_TL^p}
 \le  \|A_{-}^{(n)}\big(\psi_{+}^{(n+1)}-\psi_{+}^{(n)}\big)\|_{L^p_TL^p}
   + \|\big(A_{-}^{(n)}-A_-^{(n-1)}\big)\psi_{+}^{(n)}\|_{L^p_TL^p}
\end{equation}
 we see that the second term of this is the left hand side of 
\eqref{eqn000}.  
Since $p>1$, we choose $T>0$ small enough such that
$$
  L \equiv \Big(\frac12\Big)^{\frac{1}{p}}\left(\frac{1}2\|A_0\|_p+M^2\right)
       T^{1-\frac{1}{p}}<1,
$$
we have
\eq{\label{eqn;difference4}
 \spl{
  \|A_{-}^{(n)}(\psi_{+}^{(n+1)}-&\psi_{+}^{(n)})\|_{L^p_TL^p}\\
   \le & \frac{L}{1-L}
     \Big(\|\psi_{+}^{(n)}(A_{-}^{(n)}-A_{-}^{(n-1)})\|_{L^{p}_TL^p}
     +m\|\psi_{-}^{(n+1)}-\psi_{-}^{(n)}\|_{L^{p}_TL^p}
     \Big)
    }
}
By plugging all \eqref{eqn000}, \eqref{eqn001} and \eqref{eqn;difference4} into
\eqref{eqn;difference1}, we can replace the estimate
\eqref{eqn;difference1} by the following,
\begin{equation}\label{eqn002}
   \|\psi_{+}^{(n+1)}-\psi_{+}^{(n)}\|_{L^{\infty}_TL^p}
    \le CT^{1-\frac{1}{p}}\|\big(A_{-}^{(n)}-A_-^{(n-1)}\big) \psi_{+}^{(n)}\|_{L^p_TL^p}.
\end{equation}
We estimate by Lemma \ref{lem;bilinear},
\eq{\label{eqn;difference5}
 \spl{%
  &\|(A_{-}^{(n+1)}-A_{-}^{(n)})\psi_{+}^{(n)}\|_{L^p_TL^p}\\
   \le &  T^{1-\frac{1}{p}}
         \||\psi_+^{(n)}|^2+|\psi_-^{(n)}|^2
          -\big(|\psi_+^{(n-1)}|^2+|\psi_-^{(n-1)}|^2\big)\|_{L^p_TL^p}\\
       &\times
         \Big(\frac12\Big)^{\frac{1}{p}}
         \left(\|\psi_{+0}\|_{L^1}
               +T^{1-\frac1p}\|A_-^{(n-1)}\psi_+^{(n)}\|_{L^{p}_TL^p}
      +mT\|\psi_{-}^{(n)}\|_{L^{\infty}_TL^p} \right)
 \\
   \le & \Big(\frac12\Big)^{\frac{1}{p}}
         \left(\frac{M}2+M^2+mTM\right)
         T^{1-\frac{1}{p}}
          \||\psi_+^{(n)}|^2+|\psi_-^{(n)}|^2
          -\big(|\psi_+^{(n-1)}|^2+|\psi_-^{(n-1)}|^2\big)\|_{L^p_TL^p}. 
      }
   }
We estimate 
\eq{\label{eqn;difference7}
\spl{
   \||\psi_+^{(n)}|^2+|\psi_-^{(n)}|^2
          &-\big(|\psi_+^{(n-1)}|^2+|\psi_-^{(n-1)}|^2\big)\|_{L^p_TL^p}\\
  \le & \max\big(\|\psi_+^{(n)}\|_{L^{\infty}L^{\infty}},
                 \|\psi_+^{(n-1)}\|_{L^{\infty}L^{\infty}}\big)
        \||\psi_+^{(n)}|-|\psi_+^{(n-1)}|\|_{L^p_TL^p}\\
     & +\max\big(\|\psi_-^{(n)}\|_{L^{\infty}L^{\infty}},
                 \|\psi_-^{(n-1)}\|_{L^{\infty}L^{\infty}}\big)
        \||\psi_-^{(n)}|-|\psi_-^{(n-1)}|\|_{L^p_TL^p}\\
  \le &C\|\psi_{0\pm}\|_{L^{\infty}}
         \big(\|\psi_+^{(n)}-\psi_+^{(n-1)}\|_{L^{\infty}_TL^p}
             +\|\psi_-^{(n)}-\psi_-^{(n-1)}\|_{L^{\infty}_TL^p} \big).
 }}
Finally we obtain that 
\eq{ \label{eqn;difference10}
\spl{
  \|\psi_{+}^{(n+1)}- \psi_{+}^{(n)}\|_{L^{\infty}_TL^p}
  \le  C\left(\frac{M}2+M^2+mTM\right)
         T^{1-\frac{1}{p}}
           \big(\|\psi_+^{(n)}-\psi_+^{(n-1)}\|_{L^p_TL^p}
             +\|\psi_-^{(n)}-\psi_-^{(n-1)}\|_{L^p_TL^p}
\big).
}
}
for small $T>0$. We also obtain the  estimate similar
to $ \|\psi_{-}^{(n+1)}-\psi_{-}^{(n)}\|_{L^{\infty}_TL^p}$.

On the other hand, taking $L^{\infty}L^p$ norm for the difference of the
 integral equation
\eqn{
\spl{
 A_{\pm}^{(n+1)}-A_{\pm}^{(n)}
  = & \mp \frac12\int_0^t \big(|\psi_{\pm}^{(n)}|^2+|\psi_{\mp}^{(n)}|^2\big)ds 
      \pm\frac12\int_0^t \big(|\psi_{\pm}^{(n-1)}|^2+|\psi_{\mp}^{(n-1)}|^2\big)ds,  
    }
   }
we have from \eqref{eqn;difference7} that 
\eq{\label{eqn;difference11}
 \spl{
   \|A_{+}^{(n+1)}-A_{+}^{(n)}\|_{L^{\infty}_TL^p}
   \le& CT^{1-\frac{1}{p}}
         \big(\|\psi_+^{(n)}-\psi_+^{(n-1)}\|_{L^{\infty}_TL^p}
             +\|\psi_-^{(n)}-\psi_-^{(n-1)}\|_{L^{\infty}_TL^p} \big)^{2}.
 }
}
We have desired estimate for $A_{+}^{(n+1)}-A_{+}^{(n)}$
for sufficiently small $T>0$.
Combining two estimates \eqref{eqn;difference10} and \eqref{eqn;difference11}, 
one may easily obtain that the recurrence formula produces the uniform bound estimate 
and it is direct to conclude that the sequence 
$\{\psi_{\pm}^{(n)}, A_{\pm}^{(n)}\}_{n=1}^{\infty}$ is the Cauchy sequence 
in $L^{\infty}(0,T;L^p)\times L^{\infty}(0,T;L^p)$.
Hence there exists a pair of
the limit function $(\psi_{\pm}, A_{\pm})$ that solves the integral 
equation  \eqref{eqn;CSD-non-null}.  
This proves Proposition \ref{prop;convergence}.
\end{prf}

As we saw that $\psi_{\pm}\in L_T^{\infty}(L^p\cap L^{\infty})$ 
and $A_{\pm}\in L_T^{\infty}L^p$
imply
$|\psi_{\pm}|^2, A_{\pm}\psi_{\mp}\in L_T^{p}L^p$, 
we don't need to have extra spaces which the
 solution belongs to to obtain the uniqueness result, 
and Lipschitz continuous result 
for the solution map as well. We conclude Theorem \ref{prop3-1}.
\end{prf}

\subsection{Non-null and critical case}
Our aim in this subsection is the following,

\vskip3mm
\begin{thm}\label{prop3-2}
Let $\al=I$. 
For any 
$(\psi_{\pm0},A_{\pm0}) \in (L^1\cap L^{\infty})\times L^1$, 
there exists a unique global weak solution 
$(\psi_{\pm},A_\pm)\in C([0,\infty);L^1\cap L^{\infty})\times C([0,\infty);L^1)$ 
to \eqref{eqn;CSD-non-null}. 
The map from the initial data to the solution 
is the Lipschitz continuous from $(L^1\cap L^{\infty})
\times L^1 \to C([0,T);L^1\cap L^{\infty})
\times C([0,T);L^1)$.
\end{thm}

\begin{prf}{Theorem \ref{prop3-2}}
We use the same argument with the case of null and critical case 
in subsubsection \ref{null-critical}.
We can make $L^{1}$ norm of initial data small as much as we want
by splitting the initial data to small pieces.
We remark here that we can't make $\|\psi_{\pm0}\|_{L^{\infty}}$ 
small even if we split the support of the functions. But it doesn't 
cause a trouble. We know the solutions 
$\|\psi_{\pm}\|_{L^{\infty}_{T}L^{1}}, \|A_{\pm}\|_{L^{\infty}_{T}L^{1}}$ 
are small when the 
initial data $\|\psi_{\pm0}\|_{L^{1}}, \|A_{\pm0}\|_{L^{1}}$ are small.
So 
$\||\psi_{\pm}|^{2}\|_{L^{1}_{T}L^{1}}, \||A_{\pm}|^{2}\|_{L^{1}_{T}L^{1}}$
are small by H\"older inequality as we saw.
We follow the same ways with the proof of Proposition 
\ref{prop;convergence} but $T^{1-\frac{1}{p}}=1$. 
We can make $\left(\frac{1}2\|A_0\|_p+M^2\right)$ in 
\eqref{eqn000} and $\left(\frac{M}2+M^2+mTM\right)$
in \eqref{eqn;difference10} as small as we like. So we can obtain
for any small $C>0$, 
\begin{equation}\label{tech01}
\|\psi_{\pm}^{(n+1)}-\psi_{\pm}^{(n)}\|_{L_{T}^{\infty}L^{1}}
\le\frac{C}{2^{n}}
\end{equation}
provided that the induction assumption holds
\begin{equation*}
\|\psi_{\pm}^{(n)}-\psi_{\pm}^{(n-1)}\|_{L_{T}^{\infty}L^{1}}
+\|A_{\pm}^{(n)}-A_{\pm}^{(n-1)}\|_{L_{T}^{\infty}L^{1}}
\le\frac{1}{2^{n-1}}
\end{equation*}
and small $M$. 
From \eqref{eqn;difference11} and \eqref{tech01} with small $C$, we have
\begin{equation*}
\|A_{\pm}^{(n+1)}-A_{\pm}^{(n)}\|_{L_{T}^{\infty}L^{1}}
\le\frac{1}{2^{n+1}}
\end{equation*}
which concludes \eqref{subcritical01} with $p=1$. 
We obtain the time local existence. 

For the time global existence, we show \eqref{no-shrink2} again for 
this non-null setting. We know $\psi_{\pm}\in L^{\infty}_{T}L^{\infty}$ and
\begin{equation*}
\begin{split}
|A_{\pm}(t,x)|
&\le|A_{\pm0}(x\mp t)|+\int_{0}^{t}|\psi_{\pm}|^{2}(s,x\mp t\pm s)ds \\
&\le|A_{\pm0}(x\mp t)|+T\|\psi_{\pm}\|^{2}_{L^{\infty}_{T}L^{\infty}}.
\end{split}
\end{equation*}
These imply \eqref{no-shrink2}. We conclude the proof of 
Theorem \ref{prop3-2}.

\end{prf}

\vskip5mm\noindent
{\bf Remark.}
Since the procedure for proving the existence and 
continuity for the global solution is based on the 
local arguments, the global well-posedness for the 
\eqref{eqn;CSD} can  be generalized into 
the space of the locally uniformly class 
$L^p_{\rm loc\, unif}(\re)$, where
$$
  L^p_{\rm loc\, unif}(\re)
     =\Big\{f\in L^1_{\rm loc}(\re);
          \sup_{K\subset\subset \re} 
          \|f\|_{L^p(K)}<\infty
          \Big\}.
$$

\vskip4mm
\noindent
{\bf Acknowledgments.}  
The authors would like to express their thank to 
Professor Yoshio Tsutsumi for his valuable advice.
The work of T. Ogawa is partially supported by 
JSPS Grant-in-Aid for Scientific Research,  Basic Research
A \#20244009.


\end{document}